\documentclass[a4paper,english,10pt]{article}

\usepackage[utf8]{inputenc}
\usepackage[english]{babel}  
\usepackage{amsfonts}
\usepackage{amsmath}
\usepackage{amssymb}
\usepackage{graphicx}
\usepackage{bbold}
\usepackage{ wasysym }
\usepackage{array}
\usepackage{ textcomp }
\usepackage{hyperref}
\usepackage{hhline}
\usepackage{stmaryrd}
\usepackage{amsthm}
\usepackage[top=1in, bottom=1.25in, left=1in, right=1in]{geometry}
\usepackage{color}

\makeatletter
\def\rightharpoondownfill@{\arrowfill@\relbar\relbar\rightharpoondown}
\def\rightharpoonupfill@{\arrowfill@\relbar\relbar\rightharpoonup}
\def\leftharpoondownfill@{\arrowfill@\leftharpoondown\relbar\relbar}
\def\leftharpoonupfill@{\arrowfill@\leftharpoonup\relbar\relbar}

\newcommand{\xrightharpoondown}[2][]{\ext@arrow 0359\rightharpoondownfill@{#1}{#2}}
\newcommand{\xrightharpoonup}[2][]{\ext@arrow 0359\rightharpoonupfill@{#1}{#2}}
\newcommand{\xleftharpoondown}[2][]{\ext@arrow 3095\leftharpoondownfill@{#1}{#2}}
\newcommand{\xleftharpoonup}[2][]{\ext@arrow 3095\leftharpoonupfill@{#1}{#2}}

\newcommand{\xleftrightharpoons}[2][]{\mathrel{%
  \raise.22ex\hbox{$\ext@arrow 3095\leftharpoonupfill@{\phantom{#1}}{#2}$}%
  \setbox0=\hbox{$\ext@arrow 0359\rightharpoondownfill@{#1}{\phantom{#2}}$}%
  \kern-\wd0
  \lower.22ex\box0}}
\newcommand{\xrightleftharpoons}[2][]{\mathrel{%
  \raise.22ex\hbox{$\ext@arrow 3095\rightharpoonupfill@{\phantom{#1}}{#2}$}%
  \setbox0=\hbox{$\ext@arrow 0359\leftharpoondownfill@{#1}{\phantom{#2}}$}%
  \kern-\wd0 \lower.22ex\box0}}
\makeatother

\newcommand{\R}{\mathbb{R}}

\newcommand{\N}{\mathbb{N}}

\newcommand{\undn}{\underline{n}}
\newcommand{\beq}{\begin{equation}}
\newcommand{\eeq}{\end{equation}}

\newcommand{\bepa}{\left\{ \begin{array}{l}}
\newcommand{\eepa} {\end{array}\right.}

\newcommand{\f}{\frac}
\newcommand{\p}{\partial}

\def\restriction#1#2{\mathchoice
              {\setbox1\hbox{${\displaystyle #1}_{\scriptstyle #2}$}
              \restrictionaux{#1}{#2}}
              {\setbox1\hbox{${\textstyle #1}_{\scriptstyle #2}$}
              \restrictionaux{#1}{#2}}
              {\setbox1\hbox{${\scriptstyle #1}_{\scriptscriptstyle #2}$}
              \restrictionaux{#1}{#2}}
              {\setbox1\hbox{${\scriptscriptstyle #1}_{\scriptscriptstyle #2}$}
              \restrictionaux{#1}{#2}}}
\def\restrictionaux#1#2{{#1\,\smash{\vrule height .8\ht1 depth .85\dp1}}_{\,#2}}

\begin {document}
\author{
Martin Strugarek\thanks{AgroParisTech, Paris, France; Sorbonne Universit\'es, UPMC Univ Paris 06, CNRS, INRIA,  UMR 7598, Laboratoire Jacques-Louis Lions, \'Equipe MAMBA, 4, place Jussieu 75005, Paris, France.}
\and
Nicolas Vauchelet\thanks{Sorbonne Universit\'es, UPMC Univ Paris 06, CNRS, INRIA,  UMR 7598, Laboratoire Jacques-Louis Lions, \'Equipe MAMBA, 4, place Jussieu 75005, Paris, France.}
}
\title{Reduction to a single closed equation for 2 by 2 reaction-diffusion systems of Lotka-Volterra type.}
\date{\today}
\maketitle

\theoremstyle{plain}

\newtheorem{lemma}{Lemma}

\newtheorem{theorem}{Theorem}
\newtheorem*{nnlemma}{Lemma}
\newtheorem{proposition}{Proposition}
\newtheorem{remark}{Remark}
\newtheorem{assumption}{Assumption}
\newtheorem{corollary}{Corollary}

\begin{abstract}
We consider general models of coupled reaction-diffusion systems for interacting variants of the same species.
When the total population becomes large with intensive competition, we prove that the frequencies (\textit{i.e.} proportions) of the variants can be approached by the solution of a simpler reaction-diffusion system, through a singular limit method and a relative compactness argument.
As an example of application, we retrieve the classical bistable equation for \textit{Wolbachia}'s spread into an arthropod population from a system modeling interaction between infected and uninfected individuals.
\end{abstract}

\bigskip


\bigskip


\section{Introduction}

We are interested in modeling situations when two biological populations of the same species interact with each other, especially move, reproduce and compete.
The dynamics of these two populations are commonly described by a reaction-diffusion system of two equations in the whole space $\R^d$ ($d \geq 1$). 
In this setting, reaction terms encompass the whole interaction. Usually, they are non-linear, in order to account for competition or mutualistic interaction.
Denoting $n_1(t,x)$ and $n_2(t,x)$ the densities of each species' variant at time $t >0$ and position $x\in\R^d$,
the mathematical model reads:
\begin{equation}\label{syst1}
 \begin{cases}
  \partial_t n_1 - \nabla \cdot (A(x) \nabla n_1) &= n_1 f_1 (n_1, n_2),\\
  \partial_t n_2 - \nabla \cdot (A(x) \nabla n_2) &= n_2 f_2 (n_1, n_2),
 \end{cases}
\end{equation}
where the diffusion matrix $A$ is elliptic and the regular functions $f_1$ and $f_2$ describe
the interaction between variants.
This system is complemented with initial conditions.
Since the analysis of such systems is actually delicate, one prefers considering
the proportion of one population, for instance $p=\frac{n_1}{n_1+n_2}$.
Then the interactions are described through the dynamics of the proportion $p$
by a reaction-diffusion system:
\begin{equation}\label{syst2}
\partial_t p - \nabla \cdot (A(x) \nabla p) = p F(p).
\end{equation}

Since the pioneering works of Fisher \cite{Fisher} and Kolmogorov, Petrovskii, Piskunov \cite{KPP},
this kind of reaction-diffusion equation has been extensively studied in mathematical literature.
In particular many effort have been done to establish the existence of traveling waves and 
to describe the invasion phenomena (see e.g. \cite{Fife}, \cite{VVV}).
However, when considering systems of reaction-diffusion equations, many difficulties 
make such analysis harder. For instance, we mention the work \cite{Gardner} for competitive system.
The aim of this paper, is to focus on the link between system \eqref{syst1} and \eqref{syst2}.
More precisely, the main question we want to address is to know if solutions of system 
\eqref{syst1} can be rigorously approximated by system~\eqref{syst2} for the proportion $p$
of one species.
In our main result, we show that under suitable assumptions on the reaction terms in \eqref{syst1},
the proportion $p=\frac{n_1}{n_1+n_2}$ is close (in a sense which will be defined below) 
to a solution to system \eqref{syst2}.
More precisely, we show that when the total population becomes large with intensive competition, 
the frequency $p=\frac{n_1}{n_1+n_2}$ for system \eqref{syst1} converges to the solution 
of equation \eqref{syst2} where the non-linear function in the right hand side $F$ is explicitely
given with $f_1$ and $f_2$.
Our proof is based on a compactness argument resulting from a priori estimates.
The closest results of model reduction for competition-diffusion systems, are those of \cite{Hil.Relative} and \cite{Hil.Singular} (in bounded domains, with a specific and extensive discussion on the boundary issues).

Our first interest in this topic comes from the biological phenomenon of cytoplasmic incompatibility, caused by the endo-symbiotic bacterium \textit{Wolbachia} in some arthropod species (see \cite{Wer.Wolbachia}, \cite{BarTur.Spatial}, \cite{HugBri.Modeling}).
These bacteria have gained interest lately because of their potential use as a tool to fight arboviruses (see \cite{Hof.Successful}, \cite{Wal.wMel}).
For this situation, modeled by a reaction-diffusion system, we prove that if reaction terms scale in a proper way, then the frequency of \textit{Wolbachia} infection approaches the solution of a single closed reaction-diffusion equation, which is bistable. 
Bistable equations have been suggested long ago for this problem (see \cite{BarTur.Spatial} for an account on this topic, and \cite{Sch.Constraints} for a specific discussion). 
When these models encompass a space-dependent total population density $\rho$ (as proposed e.g. in 
\cite{Nagylaki,Barton,BarTur.Spatial}), they read
\begin{equation}
 \label{eq:BarTur}
 \partial_t p - \nabla \cdot (A(x) \nabla p) - 2 \frac{\nabla \rho}{\rho} A(x) \nabla p = p F(p).
\end{equation}
In some sense our result justifies their use thanks to a rigorous singular limit method.
We do not assume that $\rho$ and $p$ vary independently, and find that \eqref{eq:BarTur} must be corrected since $F$ is a function of $p$ and $\rho$.
We warn the reader that in order to simplify the computations, we will define a ``reduced total population density'' $n$, instead of using the total population density $\rho$ directly.



The outline of the paper is the following. In the next Section, we present the setting of the problem.
In particular the assumptions on the reaction terms and the main result are presented.
Section \ref{exwolbachia} is devoted to an example of application: the interaction between 
an infected and an uninfected mosquitoes population. A numerical illustration is also provided 
in dimension $d=1$.
The proof of our main result is provided in Section \ref{sec:proof}.
This proof relies strongly on {\it a priori} estimates that make us able to prove 
relative compactness of solutions families when a parameter describing the size 
of the population goes to $+ \infty$.
We give in Section \ref{generalization} some extension to our main result.
Finally, Section \ref{conclusion} highlights questions this work opens.

\section{Setting of the problem for typical Lotka-Volterra systems}\label{sec:model}

In this section, we first define the setting where our result applies (typical Lotka-Volterra systems), and then state it in Theorem \ref{bigthm}.

\subsection{System and assumptions}

For $\epsilon > 0$, let $f_1^{\epsilon}, f_2^{\epsilon} : \R^2 \to \R$ be two functions. 
We start from the following system in $\R^d$
\begin{equation}
\label{2pop}
 \begin{cases}
  \partial_t n_1^{\epsilon} - \nabla \cdot (A(x) \nabla n_1^{\epsilon}) &= n_1^{\epsilon} f_1^{\epsilon} (n_1^{\epsilon}, n_2^{\epsilon}),\\
  \partial_t n_2^{\epsilon} - \nabla \cdot (A(x) \nabla n_2^{\epsilon}) &= n_2^{\epsilon} f_2^{\epsilon} (n_1^{\epsilon}, n_2^{\epsilon}),
 \end{cases}
\end{equation}
with given initial data $n_i^\epsilon(t=0,x)=n_i^{\text{init}, \epsilon} \geq 0$ for $i \in \{1, 2\}$.
We assume that the matrix $A$ is elliptic and that $f_1^{\epsilon}, f_2^{\epsilon}$ are smooth enough to 
guarantee existence and uniqueness of a global solution for fixed $\epsilon>0$. More precisely,
\begin{assumption}[Ellipticity and symmetry of $A$]
 \label{assa}
 The diffusion matrix $A:\R^d\to \R^{d\times d}$ is symmetric and 
 the system \eqref{2pop} is uniformly elliptic, \textit{i.e.}
 \[
  \exists \nu_0 \in \R_+^*, \forall x, \zeta \in \R^d, \, \zeta \cdot (A(x) \zeta) \geq \nu_0 \lvert \zeta \rvert^2,
 \]
 where $\lvert \cdot \rvert$ stands for the euclidean norm in $\R^d$.
\end{assumption}

We define ``reduced total population'' $n^{\epsilon}$ and frequency (\textit{i.e.} proportion of population $1$) $p^{\epsilon}$ by
\begin{equation}
    n^{\epsilon} := \frac{1}{\epsilon} - n_1^{\epsilon} - n_2^{\epsilon}, \quad
    p^{\epsilon} := \frac{n_1^{\epsilon}}{n_1^{\epsilon}+n_2^{\epsilon}}.
    \label{def:Nnp}
\end{equation}
Since $0$ is a sub-solution for each equation in \eqref{2pop}, and since initial 
data are nonnegative, we have $n_i^\epsilon(t,\cdot)\geq 0$, for any $t\geq 0$.
By convention, we take $p^\epsilon=0$ whenever $n^\epsilon_1=n^\epsilon_2=0$.

We want to compute the limit as $\epsilon \to 0$ of the frequency $p^{\epsilon}$
under the above assumption on $A: \R^d \to \R^{d \times d}$ and on some assumptions
on the families of functions $(f_1^{\epsilon}, f_2^{\epsilon})_{\epsilon > 0}$. 

As a typical Lotka-Volterra system, we note that absence of either population of type $1$ or $2$ is a solution to this system: there is no spontaneous generation of one population from the other.
In addition, the system is positive: for non-negative initial data, $0 \leq p^{\epsilon} \leq 1$.
Now, we state our key assumptions.

\begin{assumption}[Dependence in $\epsilon$]
Functions $f_1^{\epsilon}, f_2^{\epsilon}$ are of class $\mathcal{C}^2 (\R_+^2 - \{0 \})$, and for $i \in \{1, 2\}$ there exists $F_i \in \mathcal{C}^2 (\R_+^2)$ (independent of $\epsilon > 0$) such that
\begin{equation}
\label{eq:ftoF}
f_i^{\epsilon} (n_1^{\epsilon}, n_2^{\epsilon}) = F_i (n^{\epsilon}, p^{\epsilon}).
\end{equation}
\label{asseps}
In other words, for any $n_1, n_2 \geq 0$, we may write $f_i^{\epsilon} (n_1, n_2) = F_i (\frac{1}{\epsilon} - n_1 - n_2, \frac{n_1}{n_1 + n_2})$.
\end{assumption}
From now on we drop, the superscript $\epsilon$ when it is not equivocal.

Adding the two equations in system \eqref{2pop} and using also the identity
\[
 \nabla \cdot (A(x) \nabla p) = \frac{1}{n_1 + n_2} \nabla \cdot (A(x) \nabla n_1) + \frac{p}{n_1 + n_2} \nabla \cdot (A(x) \nabla n)+ 2 \frac{1}{n_1 + n_2} \nabla p \cdot (A(x) \nabla n),
\]
we deduce, after straightforward computations, that $(n,p)$ satisfies
\begin{equation}
 \begin{cases}
\label{eq:Np}
  \partial_t n - \nabla \cdot(A(x) \nabla n) = \big(n-\f{1}{\epsilon}\big) \big( p F_1(n, p) + (1-p) F_2 (n, p) \big),  \\[2mm]
  \partial_t p - \nabla \cdot(A(x) \nabla p) + 2 \nabla p \cdot \displaystyle \frac{A(x) \nabla n}{\f{1}{\epsilon} - n}= p (1-p) \big( F_1 (n, p) - F_2(n,p) \big),
 \end{cases}
\end{equation}
complemented with well-defined initial data.
According to the first equation in \eqref{eq:Np}, it appears interesting,
when $\epsilon\to 0$, to consider the function
\begin{equation}
\label{def:H} 
H(n,p) := - p F_1(n,p) - (1-p) F_2(n,p).
\end{equation}
The following assumption guarantees existence of zeros 
$(n,p)=(h(p),p)$ for each $p\in[0,1]$ for the above function $H$.

\begin{assumption}[Nature of the interaction]
\label{assh}
In addition to Assumption \ref{asseps}, we assume
\begin{itemize}
 \item[(i)] $\exists B > 0$ such that $\forall n \geq 0, \, \forall p \in [0, 1], \, \partial_{n} H (n ,p) \leq - B$ ; 
 \item[(ii)] $\forall p \in [0, 1], \, H(0, p) > 0$.
\end{itemize}
Conditions (i) and (ii) imply that for all $p \in [0, 1]$, there exists a unique 
$n =: h(p) \in \R_+^*$ such that $H (n, p) = 0$.
We assume $H \in \mathcal{C}^2 (\R_+^2)$ 
(which is true if Assumption \ref{asseps} holds), and thus $h \in \mathcal{C}^2 (0, 1; \R)$, with $H(h(p),p) = 0$ for all $p \in [0, 1]$.
\end{assumption}

In particular, $h(0) = \f{1}{\epsilon} - \overline{n}_2^{\epsilon}$, obtained at the population 1-free equilibrium $(0, \overline{n}_2^{\epsilon})$ for \eqref{2pop}.
By Assumption \ref{assh}, this equilibrium is unique and the reduced population $n^{\epsilon}$ does not depend on $\epsilon$.

Assumption \ref{assh} may seem a little awkward, therefore we would like to point out a sufficient condition.
\begin{lemma}
 \label{ccexample}
 We assume that both $f_1^{\epsilon}$ and $f_2^{\epsilon}$ are smooth (say, of class $\mathcal{C}^2 (\R_+^2 - \{(0,0)\})$).
 We define the ``triangle'' $T_{\epsilon} = \{ (n_1, n_2) \in \R_+^2 \text{ such that } n_1 + n_2 \leq \f{1}{\epsilon} \}$.
 In addition to Assumption \ref{asseps}, if $f_1^{\epsilon}, f_2^{\epsilon}$ satisfy the following inequalities in $T_{\epsilon}$,
 \begin{equation}
 \label{Bcondition}
 \forall \undn = (n_1, n_2) \in T_{\epsilon}, \quad n_1^2 \p_{n_1} f_1^{\epsilon} (\undn) + n_1 n_2 (\p_{n_2} f_1^{\epsilon} + \p_{n_1} f_2^{\epsilon})(\undn) + n_2^2 \p_{n_2} f_2^{\epsilon} (\undn) \leq -B (n_1 + n_2)^2,
 \end{equation}
 together with a ``boundary condition'': for all $\undn \in \R_+^2$ with $\lVert \undn \rVert_1 = \f{1}{\epsilon}$,
 \begin{equation}
 \label{eq:bcexample}
 n_1 f_1^{\epsilon} (\undn) + n_2 f_2^{\epsilon} (\undn) < 0.
 \end{equation}
 Then Assumption \ref{assh} holds.
\end{lemma}
\begin{remark}
 Equation \eqref{Bcondition} on the lines $n_1= 0$ and $n_2 = 0$ means that the profiles of $f_1, f_2$ are below concave parabolic profiles.
 More generally, it ensures that the total population $n_1 + n_2$, in \eqref{eq:Np}, will start decreasing (in time) before reaching the value $\f{1}{\epsilon}$.
\end{remark}
\begin{proof}[Proof of Lemma \ref{ccexample}]
We verify each point  \textit{(i)} and  \textit{(ii)} in Assumption \ref{assh}.
 \textit{(i)} We first recall that $\p_n H = - (p \p_n F_1 + (1-p) \p_n F_2)$.
From \eqref{eq:ftoF}, we express $\p_{n_i} f_1^{\epsilon}$ and $\p_{n_i} f_2^{\epsilon}$,
 \[
    \p_{n_i} f_j^{\epsilon} = - \p_n F_j + \f{n_{3-i}}{n_1 + n_2} \p_p F_j, \qquad i,j=1,2.
 \]
Collecting these expressions yields straightforwardly
$$
  \begin{array}{ll}
\p_nH = - p \p_n F_1 - (1-p) \p_n F_2 = & \displaystyle p \f{n_1}{n_1+n_2} \p_{n_1} f_1^{\epsilon} + p \f{n_2}{n_1 + n_2} \p_{n_2} f_1^{\epsilon} \\[2mm]
&\displaystyle + (1-p) \f{n_1}{n_1 + n_2} \p_{n_1} f_2^{\epsilon} + (1-p) \f{n_2}{n_1 + n_2} \p_{n_2} f_2^{\epsilon},
  \end{array}
$$
 whence the equivalence with \eqref{Bcondition}.
 
 \textit{(ii)} For the boundary condition, we compute from \eqref{def:H}
 $$
 H (0, p) = - \big( p F_1(0, p) + (1-p) F_2(0, p) \big).
 $$
 
 Then it suffices to recall that, by definition in \eqref{eq:ftoF}, for $i \in \{ 1,2\}$, $F_i(0, p) = f_i^{\epsilon} (\f{p}{\epsilon}, \f{1-p}{\epsilon}).$
  \end{proof}

\subsection{Main result}

We are now in position to state our main result. 
We recall that we associate to any initial data $n_i^{\text{init}, \epsilon}$ the corresponding solutions of \eqref{2pop}, $(n_i^{\epsilon})$, and their relative variable $n^{\epsilon}$ and $p^{\epsilon}$, as defined in \eqref{def:Nnp}.
In addition we may define
$$
 n^{\text{init}, \epsilon} = \f{1}{\epsilon} - n_1^{\text{init}, \epsilon} - n_2^{\text{init}, \epsilon}, \quad p^{\text{init}, \epsilon} = \f{n_1^{\text{init}, \epsilon}}{n_1^{\text{init}, \epsilon}+ n_2^{\text{init}, \epsilon}}.
$$

\begin{theorem}
\label{bigthm}
We assume that Assumptions \ref{assa}, \ref{asseps}, and \ref{assh} are satisfied.
We consider the solutions of \eqref{2pop} with initial data $n_i^{\epsilon} (t=0) = n_i^{\text{init}, \epsilon} \in L^{\infty} (\R^d; \R_+)$ for $i \in \{ 1, 2 \}$.
We assume moreover that there exists $p^{\text{init}} \in L^2 (\R^d)$ such that
\begin{equation}\label{hyp:init}
p^{\text{init},\epsilon} \underset{\epsilon \to 0}{\rightharpoonup} p^{\text{init}} \mbox{ in } L^2 (\R^d)-\mbox{weak}, 
\quad  n^{\text{init},\epsilon} - h(0) \in L^2\cap L^\infty (\R^d),
\end{equation}
with uniform bounds in $\epsilon > 0$.

Then, for all $T > 0$, defining $\mathcal{H}^1_T = L^2 (0, T ; L^2 (\R^d))$ and $\mathcal{H}^2_T = L^2 (0, T ; H^1 (\R^d))$, we have the convergence 
\begin{equation}\label{conv}
 \begin{cases}
    p^{\epsilon} \xrightarrow[\epsilon \to 0]{} p^0 \text{ strongly in $\mathcal{H}^1_T$,} \text{ weakly in } \mathcal{H}^2_T,\\
    n^{\epsilon} - h(p^\epsilon) \xrightarrow[\epsilon \to 0]{} 0 \text{ strongly in $\mathcal{H}_T^1$,} \text{ weakly in } \mathcal{H}^2_T,
 \end{cases}
\end{equation}
where $p^0$ is the unique solution of the following initial value problem
\begin{equation}\label{limeq}
 \begin{cases}
  \partial_t p^0 - \nabla \cdot (A(x) \nabla p^0) = p^0 F_1 (h(p^0), p^0),\\
  p^0 (t=0) = p^{\text{init}}.
 \end{cases}
\end{equation}
\end{theorem}

This result asserts that, locally in time, the proportion of the first population, $p$,
solution to system~\eqref{2pop}, under suitable assumption 
on the reaction term and on the initial data, is close to the solution
of a single reaction-diffusion system \eqref{limeq}.
This latter system have been intensively studied, in particular existence of 
traveling waves, describing propagation phenomena (see e.g. \cite{Fife}, \cite{VVV}).
The main interest in this reduction process is that since the behavior of solutions to the scalar equation \eqref{limeq}
is well-known.
Therefore we can deduce, for small values of $\epsilon$, 
the local in time behavior of solutions to \eqref{2pop}.

We observe that the limit reaction term $r(p) := p F_1 (h(p), p)$ in \eqref{limeq} 
satisfies
$$r(0)=0, \qquad r(1) = F_1 (h(1), 1) = 0,$$
 because $H(h(p), p ) = 0 = - p F_1(h(p),p) - (1-p) F_2(h(p),p)$.
 It means that the states $p^0=0$ (only population $2$) and $p^0=1$ 
 (only population $1$) are equilibria for this system.

 Moreover, 
 $$
 r'(0) = F_1(h(0), 0)
 $$
 and
 $$
 r'(1) = h'(1) \p_{n} F_1(h(1), 1) + \p_p F_1 (h(1), 1).
 $$
 Hence under some direct sign assumptions on $F_1$ and $\p_p F_1$, 
 the equilibria $0$ and $1$ for $p$ can be made stable in the limit equation, if $r'(0)$ and $r'(1)$ are negative.
 In particular, in the example in Section \ref{exwolbachia}, the function $r$
 is bistable.

\begin{remark}
\label{rmkinit}
The assumption $n^{\text{init},\epsilon} - h(0)$ uniformly bounded with respect to $\epsilon$ in $L^2 (\R^d)$ together with the uniform bound of $p^{\text{init}, \epsilon}$ in $L^2(\R^d)$
imply, thanks to Assumption \ref{assh}, that $ n^{\text{init}} - h (p^{\text{init}, \epsilon})$ is bounded
in $L^2 (\R^d)$, uniformly in $\epsilon > 0$.
Indeed, $h$ is Lipschitz on $(0,1)$ and by the triangle inequality, we have
$\|n^{\text{init}} - h (p^{\text{init}, \epsilon})\|_{L^2}\leq 
\|n^{\text{init},\epsilon} - h(0)\|_{L^2} + \|h\|_{Lip} \|p^{\text{init}, \epsilon}\|_{L^2}$.
 \end{remark}

\begin{remark}
\label{compactsupport}
One might be interested by the effect of the local introduction of a variant into a population at 
equilibrium. In this situation, at the time of introduction, variant $2$ is at equilibrium whereas 
the introduction of variant $1$ is modeled by a compactly supported continuous nonnegative function 
$\phi$.
Then we have $n_2^{\text{init},\epsilon} = \frac 1\epsilon -h(0)$ on $\R^d\setminus \mbox{supp } \phi$,
and we set $n_1^{\text{init},\epsilon} = \phi n_2^{\text{init},\epsilon}$.
Then, $p^{\text{init}}=\frac{\phi}{1+\phi}$ and assumption~\eqref{hyp:init}
in Theorem \ref{bigthm} boils down to assuming that
$$
\frac 1\epsilon -(1+\phi) n_2^{\text{init},\epsilon}-h(0) \mbox{ is uniformly bounded with respect to } 
\epsilon \mbox{ in } L^\infty(\R^d).
$$
 \end{remark}

Finally, we mention that we can relax the assumption \eqref{hyp:init} by assuming that
the sequence $(p^{\text{init},\epsilon})_\epsilon$ is uniformly bounded with respect to $\epsilon$
in $L^2(\R^d)$ instead of assuming its convergence. In fact, we can extract a subsequence
of $(p^{\text{init},\epsilon})_\epsilon$ that converges weakly towards $p^{\text{init}}$ 
and the result applies. But the uniqueness of the weak limit $p^{\text{init}}$ is not guaranteed
and therefore, the result in Theorem \ref{bigthm} is available only up to an extraction 
of a subsequence.

\section{Application to a biological example}
\label{exwolbachia}

\subsection{Presentation of the model}
We consider the case of \textit{Wolbachia} in arthropod species (for the biology of this bacterium, see \cite{Wer.Wolbachia} ; for mathematical modeling, see \cite{BarTur.Spatial}, \cite{Fen.Solving}, \cite{HugBri.Modeling}, \cite{ChaKim.Modeling}). 
It is an endo-symbiont that is maternally transmitted, causes cytoplasmic incompatibility (CI), and has several other effects on its host.
Here, we understand CI as a mechanism through which one of the possible crossings is less viable. More precisely, if an uninfected female is fertilized by an infected male, a fraction only of its eggs will eventually hatch and give birth to viable larvae. 
For more details about CI, we refer to \cite{Wer.Wolbachia}.
In the case of \textit{Aedes} mosquitoes, \textit{Wolbachia} reduces lifespan, changes fecundity and blocks the development of dengue virus (see  \cite{Moreira}, \cite{Wal.wMel}, \cite{Joa.Distribution}).
It is then a potential biological tool to fight dengue epidemics.
However, it does not change the way mosquitoes move. 
Therefore, in order to model a \textit{Wolbachia} invasion (assessed in the field in \cite{Hof.Successful}) we are precisely in our setting. Several (two) \emph{variants of the same species} interact with each other in a complex way.

Specifically, we define the uninfected death rate $d_u$. This rate is multiplied by $\delta > 1$ for infected mosquitoes: $d_i = \delta d_u$.
We also define an uninfected fecundity $F_u$ for uninfected mosquitoes, $F_i=(1-s_f)F_u$ for infected mosquitoes ; a resource parameter $\sigma$ ;
and a CI parameter $0 < s_h \leq 1$, which means that a fraction $s_h$ of uninfected females' eggs fertilized by infected males won't hatch. Parameters $\delta$, $s_f$ and $s_h$ have been estimated in several cases and can be found in the literature (see \cite{BarTur.Spatial} and references therein).
We will always assume $s_h>s_f$. (In practice, we usually have $s_f$ close to $0$ and $s_h$ close to $1$).
Let us denote $n_i(t,x)$, resp. $n_u(t,x)$, the density of the infected, resp. uninfected,
mosquitoes at time $t\geq 0$, position $x\in\R^d$.

Several models have been written, using these parameters. 
In \cite{ChaKim.Modeling} (if we ignore the drift speed $v \in \R^d$ they used, which amounts at a change of coordinates) one find
\begin{equation}
 \label{eq:ChaKim}
 \begin{cases}
    &\p_t n_i - \nabla \cdot (A(x) \nabla n_i) = n_i (1 - \sigma(n_u + n_i)) - d_u n_i, \\
    &\p_t n_u - \nabla \cdot (A(x) \nabla n_u) = n_u F_u (1-s_h \f{n_i}{n_u + n_i})  (1 - \sigma(n_u + n_i)) - d_u n_u.
 \end{cases}
\end{equation}
In this model, $\delta = 1$ and variables are scaled so that $F_u (1-s_f) = F_i = 1$. Here the reduced population is defined by $n = \f{1}{\sigma} - (n_i + n_u)$.
The corresponding dynamics in $(n, p)$ for \eqref{eq:ChaKim} is written
\begin{equation}
\label{eq:ChaKimNp}
 \begin{cases}
    &\p_t n - \nabla \cdot (A(x) \nabla n) =  \big( \sigma n(p + F_u (1-p) (1-s_h p)) - d_u\big),\\
    &\p_t p - \nabla \cdot (A(x) \nabla p) + 2 \frac{\nabla n}{n} A(x) \nabla p = \sigma n p (1-p) \big(d_u - F_u (1-s_hp)\big),
 \end{cases}
\end{equation}
In \eqref{eq:ChaKimNp}, the reaction term for $p$ depends on $n$ merely for its intensity (it is a multiplicative factor).
In particular, the unstable steady state (defining, in some sense, \textit{a} possible ``threshold for invasion'') is equal to $ \frac{1}{s_h} (1 - \frac{d_u}{F_u})$ does not depend on $n$.

To further reduce this class of models and prove the convergence towards \eqref{eq:BarTur}, we introduce the parameter $\epsilon$ to characterize the high fertility and competition that result in a carrying capacity of order $\f{1}{\epsilon}$.
Then we propose the following generalization of \eqref{eq:ChaKim}, which incorporates also
the different death rate and the reduction of fecundity,
\begin{equation}
\label{wolsys1}
 \begin{cases}
    \p_t n_i - \nabla \cdot (A(x) \nabla n_i) &= (1-s_f) F_u n_i \big(\f{1}{\epsilon} - \sigma (n_i + n_u) \big) - \delta d_u n_i, \\
    \p_t n_u - \nabla \cdot (A(x) \nabla n_u) &= F_u n_u (1 - s_h \f{n_i}{n_i + n_u} ) \big(\f{1}{\epsilon} - \sigma (n_i + n_u) \big) - d_u n_u, \\
 \end{cases}
\end{equation}

Straightforwardly, we can compute the equilibria for the associated dynamical system.
\begin{lemma}
As soon as $s_f + \delta - 1 < \delta s_h$, there are four distinct equilibria associated with \eqref{wolsys1} in the non-negative quadrant.
\begin{itemize}
\item \textit{Wolbachia} invasion steady state
$(n^*_{iW}, n^*_{uW}) := (\f{1}{\sigma \epsilon}-\f{d_u}{F_u} \f{\delta}{1 - s_f},0)$
is stable;
\item \textit{Wolbachia} extinction steady state 
$(n^*_{iE}, n^*_{uE}) := (0,\f{1}{\sigma \epsilon}-\f{d_u}{F_u})$ 
is stable;
\item The co-existence steady state 
$(n^*_{iC}, n^*_{uC}) := \big(\big(\f{1}{\sigma\epsilon}-\f{d_u}{F_u} \f{\delta}{1-s_f}\big)
\f{\delta-(1 - s_f)}{\delta s_h},
\big(\f{1}{\sigma\epsilon}-\f{d_u}{F_u} \f{\delta}{1-s_f}\big)
\f{\delta(s_h-1)+(1 - s_f)}{\delta s_h} \big)$ 
is unstable;
\item The steady state $(0,0)$ is unstable.
\end{itemize}
\label{equilibria}
\end{lemma}

\subsection{Large population asymptotic}
We perform the limit $\epsilon\to 0$ for system \eqref{wolsys1}.
To recover notations from Theorem \ref{bigthm}, we identify $n_1 = n_i$, $n_2 = n_u$.
As above we define the reduced quantity $n=\f{1}{\sigma\epsilon}-(n_1+n_2)$ and
$p=\f{n_1}{n_1+n_2}$.
Then with the notations in Section \ref{sec:model}, one has
\begin{align*}
F_1 (n, p) &= \sigma n (1-s_f) F_u - \delta d_u,\\
F_2 (n, p) &= \sigma n F_u (1 - s_h p) - d_u.
\end{align*}
Therefore, by definition \eqref{def:H}, we compute
\begin{align*}
 H(n, p) &= -p ( \sigma n(1-s_f) F_u - \delta d_u) -(1-p) (\sigma n F_u (1 - s_h p) -d_u)\\
 &= - \sigma F_u n (s_h p^2 - (s_f + s_h) p + 1) + d_u ((\delta - 1) p + 1).
\end{align*}
And Assumption \ref{asseps} is satisfied. Then, Assumption \ref{assh} is easy to check since 
$H (0, p) = d_u ((\delta - 1) p + 1) > 0$ and using the fact that the polynomial 
$x\mapsto s_h x^2-(s_f+s_h)x +1$ is minimal for $x=\frac{s_f+s_h}{2s_h}$, we have
\begin{equation*}
\p_{n} H (n, p) = - \sigma F_u (s_h p^2 - (s_f + s_h) p + 1) \leq - \sigma F_u (1 - \f{(s_f + s_h)^2}{4 s_h}) < 0,
\end{equation*}
since we have $s_f<s_h$.
We notice also that this computation implies that the above second order polynomial
in $p$ is always away from $0$ on $[0,1]$.
Moreover, recalling the definition $H(n,p)=0$ if and only if $n=h(p)$ from
Assumption \ref{assh}, we can compute $ h(p) = \f{d_u}{\sigma F_u} \f{(\delta - 1)p + 1}{s_h p^2 - (s_f + s_h) p + 1}.$
Under Assumption \ref{assa} on $A$, Theorem \ref{bigthm} applies, and $p^{\epsilon}$ converges towards the solution of the following equation
\begin{equation}
\left\{
 \begin{array}{rcl}
    \p_t p^0 - \nabla \cdot (A(x) \nabla p^0) &= &r(p^0), 
    \\
    p^0 (t=0) &= &p^{\text{init}},
 \end{array}
 \right.
 \label{limeqWolba}
\end{equation}
and the reaction term writes
$$
r(p) = \delta d_u s_h \f{p (1 - p)(p-\theta)}{s_h p^2 - (s_f + s_h) p + 1}, 
\qquad \theta = \f{s_f + \delta - 1}{\delta s_h},
$$
which is bistable provided $\delta$ satisfies the condition from Lemma \ref{equilibria}:
\beq
s_f+\delta -1<\delta s_h.
\label{conddelta}
\eeq
If $\delta = 1$, we find the ubiquitous value $\theta (= p^*) = \frac{s_f}{s_h}$,
which corresponds to the model of spacial spread of Wolbachia proposed in \cite{BarTur.Spatial}.
In addition, this expression is coherent with the one in \cite{Sch.Constraints} for general $\delta$.
Even though the equation for $p$ has already been suggested for a while, as far as we know, no convergence result as ours had been proved before from a two-populations model to the bistable equation.

A direct application of Theorem \ref{bigthm} establishes that, in the limit $\epsilon \to 0$, 
the derivation in \cite{BarTur.Spatial} holds true in a strong topology.
\begin{corollary}
Assume that $A$ satisfies Assumption \ref{assa}.
Given $n_1^{\text{init},\epsilon}$ and 
$n_2^{\text{init},\epsilon}$ such that there exists $p^{\text{init}}\in L^2(\R^d)$ such that 
$p^{\text{init},\epsilon}\rightharpoonup p^{\text{init}}$ as $\epsilon\to 0$ in $L^2(\R^d)$-weak 
and $\f{1}{\sigma\epsilon} - n_1^{\text{init}, \epsilon} - n_2^{\text{init}, \epsilon} - \f{d_u}{\sigma F_u} \in L^2 \cap L^{\infty} (\R^d)$ with uniform bounds in $\epsilon > 0$,
then Theorem \ref{bigthm} applies and the solutions $(n_i^{\epsilon}, n_u^{\epsilon})_{\epsilon > 0}$ 
of \eqref{wolsys1}
  satisfy the convergence result in \eqref{conv} where the limiting equation is given in \eqref{limeqWolba}.
  \label{cor1}
\end{corollary}

\subsection{Numerical illustration}
\label{sec:num}

\begin{figure}[h!]
 \includegraphics[width=.5\linewidth]{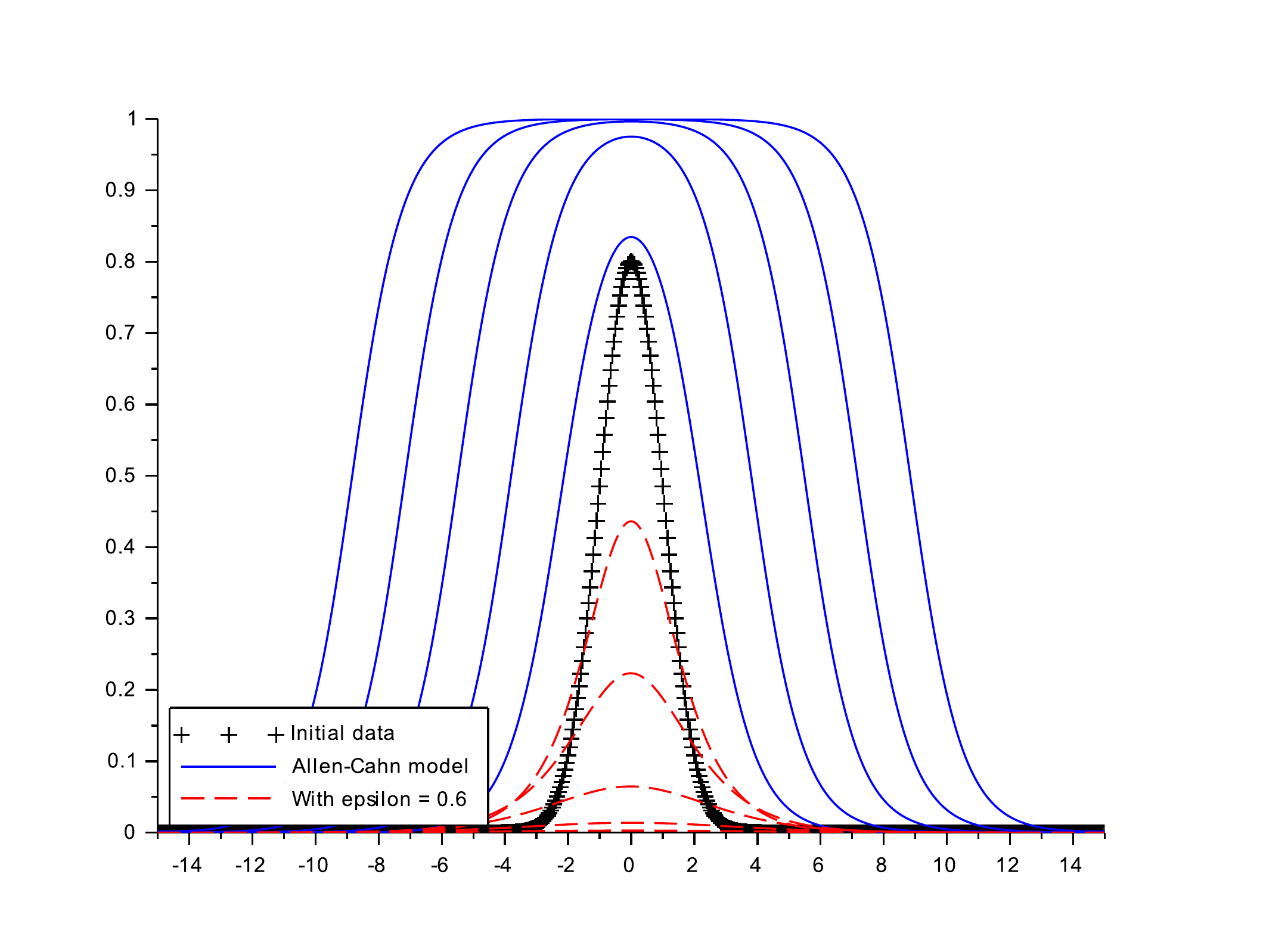}
 \includegraphics[width=.5\linewidth]{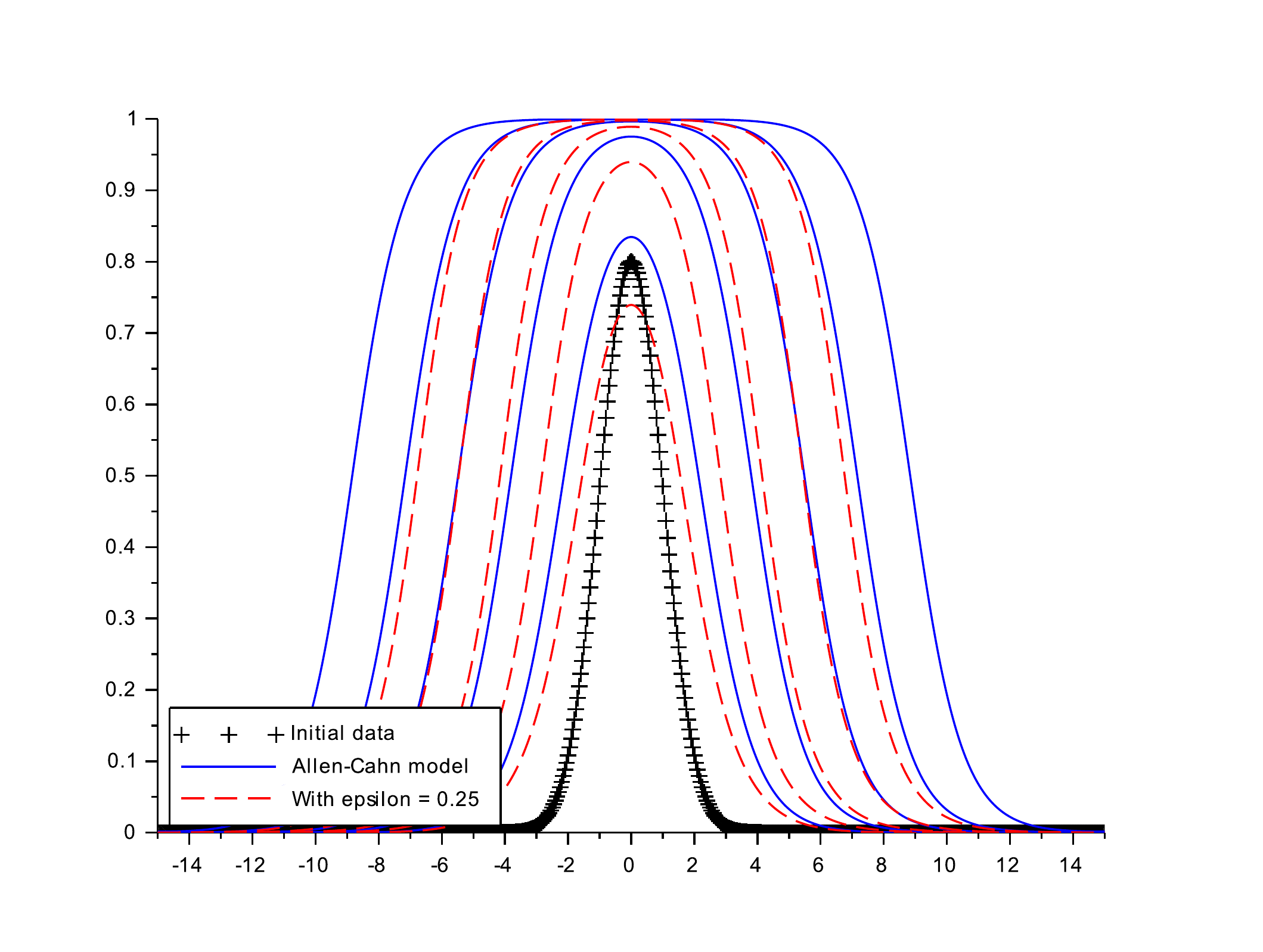} 
 \includegraphics[width=.5\linewidth]{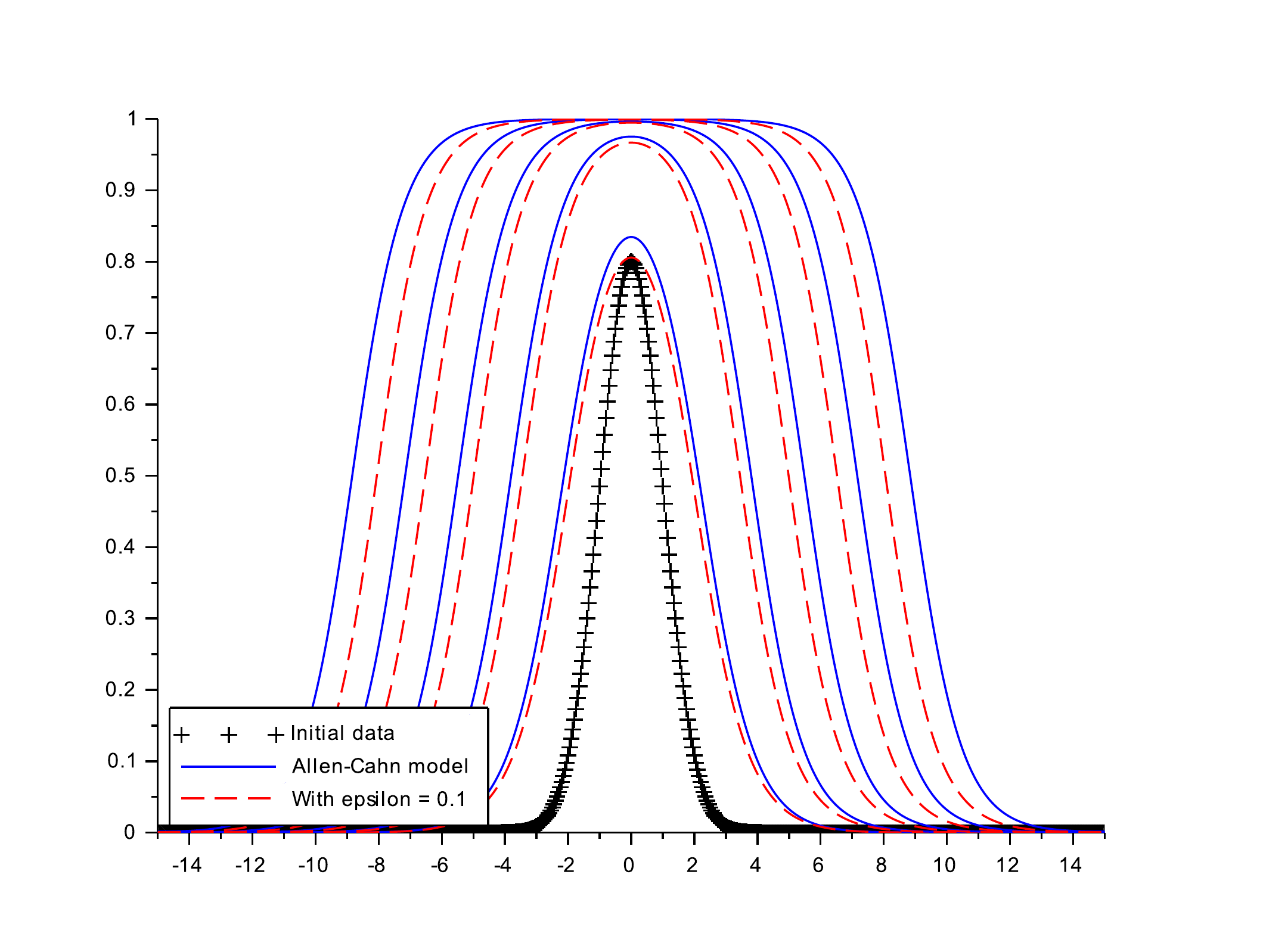}
 \includegraphics[width=.5\linewidth]{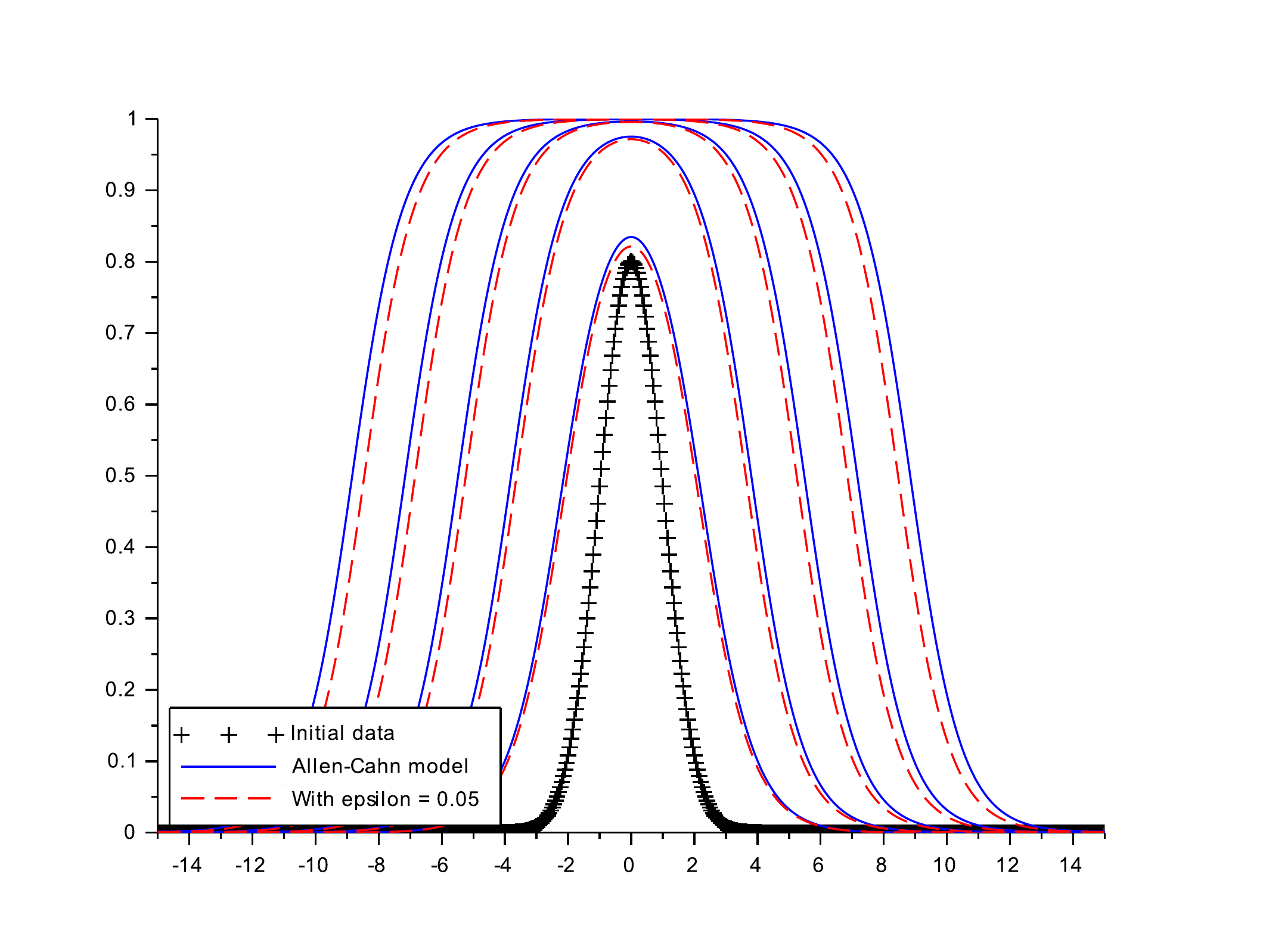}
 \caption{Initial data (+) creating a traveling wave in the limit system (blue) and convergence of the two-population solution (dashed red) as $\epsilon$ diminishes.}
 \label{fig:convergence}
\end{figure}

A numerical illustration of this convergence result is shown in Figure \ref{fig:convergence}. Parameters are fixed according to biologically relevant data (freely adapted from \cite{Foc.Dynamic}). Time unit is the day, and parameters per day are $F_i = F_u = 1.12$ (hence $s_f=0$), $d_u = .27$ and $d_i = .3$, then $\delta = \frac{d_i}{d_u}=\frac{10}{9}$.
We choose $s_f=.1$ and $s_h = .8$.
We take $\sigma=1$, and $A(x) \equiv .1$, which amounts at choosing a space scale.

We discretize the one-dimensional computational domain $[-15;15]$ 
with space step $\Delta x = .05$ and take a time step $\Delta t = .005$.
The reaction diffusion equations are discretized thanks to semi-implicit finite 
difference scheme, the diffusion operator being treated implicitly 
(to avoid too restrictive stability conditions),
while the reaction term is treated explicitly.
Curves are plotted every $5000$ iterations, at times (in days) 
$T_1=25$, $T_2=50$, $T_3 = 75$, $T_4 = 100$ and $T_5 = 125$.
We display 4 numerical tests with the same initial data $p^{\text{init}}$ 
compactly supported, plotted in pluses (+).
The blue lines represent the solution of the limiting system \eqref{limeqWolba}.
In dashed red lines are plotted the computed fraction $p=\frac{n_i}{n_i+n_u}$ 
where $(n_i,n_u)$ solves the system of two populations \eqref{wolsys1}.
In this computation the initial data are taken as noticed in Remark \ref{compactsupport}:
we consider the compactly supported function $\phi=p^{\text{init}}/(1-p^{\text{init}})$
and for given $n_u^{\text{init}}$ we take $n_i^{\text{init}}=\phi n_u^{\text{init}}$.
We observe that the solution of the limiting bistable system \eqref{limeqWolba}
exhibits a traveling front which propagates into the whole domain. 
Then the numerical results for 4 different values of the parameter $\epsilon$ 
are represented.
For large populations, we observe that as $\epsilon$ goes to $0$ 
(recall that the order of magnitude of the population size is $\f{1}{\sigma \epsilon}$), 
the solution to the whole system \eqref{wolsys1} gets closer to the one of the limiting system.
However, for $\epsilon=.6$, the introduced population goes extinct, and $p$ does not behave as in the limiting model. This illustrates how the 2 by 2 system qualitatively differs from the limit reaction-diffusion equation.

An additional conclusion we can draw from Figure \ref{fig:convergence} is that our approximation result will always be \textit{local in time}.
Indeed, for small $\epsilon$ we see a traveling wave appear in dashed red, that has a \textit{slower} speed than the blue one. 
Hence the norm of their difference will be constantly growing in time.

\section{Proof of convergence}\label{sec:proof}

This Section is devoted to the proof of Theorem \ref{bigthm}.
We write the system of equations satisfied by $(n^\epsilon, p^\epsilon)$
\begin{equation}
\label{N1p}
 \begin{cases}
  \partial_t n^\epsilon - \nabla \cdot (A(x) \nabla n^\epsilon) = \big(\f{1}{\epsilon} - n^\epsilon \big)H (n^\epsilon, p^\epsilon), \\[2mm]
  \partial_t p^\epsilon - \nabla \cdot (A(x) \nabla p^\epsilon) + 2 \epsilon \nabla p^\epsilon \cdot \frac{A(x) \nabla n^\epsilon}{1 - \epsilon n^\epsilon} = p^\epsilon (1-p^\epsilon) \big(F_1 - F_2 \big)(n^\epsilon, p^\epsilon),\\[2mm]
  n^\epsilon(t=0) = n^{\text{init},\epsilon},  \qquad   p^\epsilon (t=0) = p^{\text{init},\epsilon} \,.
 \end{cases}
\end{equation}
We recall that the initial data are assumed to satisfy \eqref{hyp:init}.
Then the sequence $(p^{\text{init},\epsilon})_\epsilon$ is bounded uniformly in $\epsilon$ 
in $L^2(\R^d)$ and $(n^{\text{init},\epsilon}-h(0))_\epsilon$ is bounded uniformly in $\epsilon$ 
in $L^2\cap L^\infty(\R^d)$.
The proof of Theorem \ref{bigthm} relies strongly on a sequence of a priori estimates
uniform in $\epsilon$, which give compactness and allow to pass to the limit in the 
equation for $p^\epsilon$.

From now on, we will drop the superscript $\epsilon$ in the notations.

\subsection{Estimates}

For $\epsilon>0$ fixed, existence of solutions to \eqref{N1p} is classical (see e.g. \cite{Perthame}).
Now we establish some a priori estimates uniform in $\epsilon>0$.
First, we have the following $L^\infty$ bounds.
\begin{lemma}
\label{infbounds}
 Under the assumptions of Theorem \ref{bigthm}, for any positive initial data, the unique solution $(p, n)$ to \eqref{N1p} satisfies
 $$
 \forall t> 0, x \in \R^d, \, 0 \leq p(t, x) \leq 1
 $$
 and $n \in L^{\infty}(\R_+ \times \R^d)$. Moreover, there exists $\epsilon_0 > 0$ such that the $L^{\infty}$ bound on $n$ is uniform in $\epsilon_0 > \epsilon > 0$.
\end{lemma}

\begin{proof}
 As stated before, positivity of $n_1, n_2$ is straightforward and implies the uniform
bounds on $p$ in $L^\infty$.
 
Using Stampacchia’s method for the bound on $n$, we notice that, from Assumption \ref{assh}, for all $p\in [0,1]$, $\big(\f{1}{\epsilon} - n \big)H (n, p)$ is positive for $n$ between $0$ and $h(p)$ and negative afterwards until $\f{1}{\epsilon}$.
 Then, for $\tilde{K} = \max_{p \in [0, 1]} h(p)$, we define $y(t) = \int_{\R^d} \big( n (t, x) - \tilde{K} \big)_+ dx$. Multiplying the equation on $n - \tilde{K}$ by $1_{n > \tilde{K}}$ and
 integrating over $\R^d$ gives, for $\epsilon < \f{1}{\tilde{K}}$
 $$
 \f{d}{dt} y(t) + \int_{\R^d} \nabla (n - \tilde{K})_+ \cdot A(x) \nabla (n - \tilde{K})_+ dx \leq 0.
 $$
 And in particular, $\f{d}{dt}{y} < 0$. 
 Since, from Assumption \eqref{hyp:init}, $n^{\text{init}}$ is bounded in $L^\infty$ uniformly
with respect to $\epsilon$, we can pick $\tilde{K}$ such that $\tilde{K} > \lVert n^{\text{init}} \rVert_{\infty}$.
 Then $y(0) = 0$. We deduce that $y \equiv 0$.
 
 To conclude, the result is proved with
 $$
 \epsilon_0 = \big(\max \, ( \max_{p \in [0, 1]} h(p), \lVert n^{\text{init}} \rVert_{\infty})  \big)^{-1}.
 $$
 \end{proof}

Now, we aim at getting the following boundedness result.
\begin{proposition}
\label{lem:universalEstimates}
Let $T>0$. Under the assumptions in Theorem \ref{bigthm}, we define $M := n - h(p)$.
Then, there exists $\epsilon_0 > 0$ such that $M$ and $p$ are uniformly bounded in $\mathcal{H}^1_T \cap \mathcal{H}^2_T$, for all $\epsilon \leq \epsilon_0$.
\end{proposition}

We recall that the function $h$ is defined in Assumption \ref{assh} and belongs to $\mathcal{C}^2([0,1])$.
Then we may define
\begin{equation}\label{def:h0}
h_0 = \|h\|_{L^\infty([0,1])}, \quad h'_0 = \|h'\|_{L^\infty([0,1])}, \quad 
h''_0 = \|h''\|_{L^\infty([0,1])}.
\end{equation}
We notice that, by definition and from Lemma \ref{infbounds}, we have that $M$ is uniformly bounded 
in $L^\infty$ for $\epsilon\leq \epsilon_0$.
The proof of this result relies on estimates on $p$ and $M$, and we postpone the proof of Proposition \ref{lem:universalEstimates} after proving them in the two following technical Lemma. The first one is for $p$.
\begin{lemma}
\label{lem:pestimate}
There is a positive constant $K$ independent of $\epsilon$ such that $\forall \epsilon > 0$,
\[
	\f{1}{2}\f{d}{dt} \int_{\R^d} p^2 dx + (1- \epsilon C_1) \int_{\R^d} \nabla p A(x) \nabla p dx \leq \epsilon C_2 \int_{\R^d} \nabla M A(x) \nabla M dx + K \int_{\R^d} p^2 dx,
\]
where $C_1 = 2 \big( 1 + \f{h'_0}{2} + (h'_0)^2 \big)$ and $C_2 = 2(1 + \f{h'_0}{2})$.
\end{lemma}
\begin{proof}
We multiply by $p$ the equation satisfied by $p$ in \eqref{N1p}, and integrate over $\R^d$
\begin{multline}\label{estimp1}
\f{1}{2} \f{d}{dt} \int_{\R^d} p^2 dx+ \int_{\R^d} \nabla p A(x) \nabla p dx + 2 \epsilon \int_{\R^d} \f{p}{1-\epsilon n} \nabla p \cdot A(x) \nabla n dx 
\\ \leq \int_{\R^d} p^2 (1-p) \big(F_1 - F_2 \big) (n, p) dx.
\end{multline}
Thanks to Lemma \ref{infbounds}, we know that $\f{p}{1 - \epsilon n}$ is well-defined for $\epsilon$ small enough, and the denominator is uniformly positive.
Hence we may use a Cauchy-Schwarz inequality,
\[
\int_{\R^d} \f{p}{1 - \epsilon n} \nabla n A(x) \nabla p dx \leq 
\f{1}{2} \int_{\R^d} \f{p}{1- \epsilon n} \nabla p A(x) \nabla p dx 
+ \f{1}{2} \int_{\R^d} \f{p}{1-\epsilon n} \nabla n A(x) \nabla n dx.
\]
Since $n = M + h(p)$, we have $\nabla n = \nabla M + h'(p)\nabla p$.
We may also write,
\[
\int_{\R^d} \nabla n A(x) \nabla n dx \leq ((h'_0)^2 + \f{h'_0}{2})\int_{\R^d} \nabla p A(x) \nabla p dx + (1 + \f{h'_0}{2})\int_{\R^d} \nabla M A(x) \nabla M dx.
\]
Now, collecting these inequalities for $\epsilon$ small enough, such that $\f{p}{1 - \epsilon n} \leq 2$ yields
\begin{multline*}
\f{1}{2} \f{d}{dt} \int_{\R^d} p^2 dx + \Big(1 - 2 \epsilon \big( 1+ h'_0(h'_0 + \f{1}{2})\big) \Big) \int_{\R^d} \nabla p A(x) \nabla p 
\\
\leq 2 \epsilon \big( 1 + \f{h'_0}{2}\big) \int_{\R^d} \nabla M A(x) \nabla M + K_F \int_{\R^d} p^2 dx,
\end{multline*}
where 
\begin{equation}\label{def:KF}
K_F := \sup \{  \lvert (1 - p) (F_1 - F_2) (n, p) \rvert, \text{ as }\lvert n \rvert \leqslant \sup_{0<\epsilon\leqslant\epsilon_0} \|n\|_{L^\infty} \text{ and } 0 \leqslant p \leqslant 1\}.
\end{equation}
Thanks to Lemma \ref{infbounds} and the continuity of the functions 
$F_1$ and $F_2$, the constant $K_F$ is finite.
This is the expected estimate.
\end{proof}

Similarly, on $M:=n-h(p)$,
\begin{lemma}
\label{lem:Mestimate}
There are positive constants $C_3, C_4, C_5, C_6$ independent of $\epsilon$ such that, for all $\epsilon > 0$,
\begin{multline*}
	\f{1}{2} \f{d}{dt} \int_{\R^d} M^2 dx + (1 - \epsilon C_3) \int_{\R^d} \nabla M A(x) \nabla M dx 
\\ \leq \big(C_4- \f{C_5}{\epsilon} \big) \int_{\R^d} M^2 dx + C_6 \big(\int_{\R^d} p^2 dx + \int_{\R^d} \nabla p A(x) \nabla p dx \big),
\end{multline*}
where $C_3 = \f{h'_0}{2}$, $C_4 = B(h_0 + \tilde{K}') + \f{h'_0 K_F}{2}$, $C_5 = B$ and $C_6 = \max\big(\f{h'_0 K_{F}}{2}, \tilde{K}' h''_0 + \epsilon h'_0(h'_0 + \f{1}{2}) \big)$, and $\tilde{K}', K_{F}$ are positive constants defined in \eqref{def:KF} and \eqref{def:K'}.
\end{lemma}
\begin{proof}
The quantity $M$ satisfies the following equation (obtained from \eqref{N1p})
\begin{multline}\label{eqM1}
  \p_t M - \nabla \cdot (A(x) \nabla M) = \p_t n - \nabla \cdot (A(x) \nabla n) - h'(p) \p_t p + \nabla \cdot (h'(p) A(x) \nabla p)
  \\= \big( \f{1}{\epsilon} - M - h(p) \big)H (M+h(p),p) + h''(p) \nabla p \cdot ( A(x) \nabla p )
  \\ - h'(p) \Big( (1-p) (- F_2 - H) (M + h(p), p) 
  \\- 2 \epsilon\f{1}{1 - \epsilon n} (A(x) \nabla M + h'(p) A(x) \nabla p) \cdot \nabla p \Big),
 \end{multline}
it is associated with an initial data $M^{\text{init}} = n^{\text{init}} - h(p^{\text{init}})$ bounded in $L^2 (\R^d)$.
Indeed, as noted in Remark \ref{compactsupport}, we have,
\[
\lvert n^{\text{init}} - h(p^{\text{init}}) \rvert \leq \lvert n^{\text{init}} - h(0) \rvert + \lvert h(0) - h(p^{\text{init}})| \leq \lvert n^{\text{init}} - h(0) \rvert + h'_0 \lvert p^{\text{init}} \rvert.
\] 
Moreover, from \eqref{hyp:init}, $p^{\text{init}}$ is bounded in $L^2 (\R^d)$ and $n^{\text{init}} - h(0)$ is bounded in $L^2(\R^d)$, with uniform bounds in $\epsilon$. It implies the uniform bound of $M^{\text{init}}$ in $L^2(\R^d)$.

Now, we assume that $\epsilon$ is small enough, so that the term $\f{1}{\epsilon} - M - h(p)$ remains positive (this is possible thanks to Lemma \ref{infbounds}).
We multiply by $M$ equation \eqref{eqM1}, and integrate over $\R^d$
\begin{multline}\label{eqM2}
 \f{d}{dt} \int_{\R^d} M^2 dx + \int_{\R^d} \nabla M \cdot ( A(x) \nabla M) dx 
\\ \leq - B \int_{\R^d} M^2 \big( \f{1}{\epsilon} - M - h(p) \big) dx
 \\ + \int_{\R^d} 2 M h'(p) \epsilon\f{1}{1 - \epsilon n} (A(x) \nabla M + h'(p) A(x) \nabla p) \cdot \nabla p  dx
 \\ + \int_{\R^d} M h''(p) \nabla p \cdot (A(x) \nabla p) dx
 \\ - \int_{\R^d} M h'(p) (1-p) (- F_2 - H) (M+ h(p), p) dx
\end{multline}
Since $\p_{n} H \leq - B$ (Assumption \ref{assh}), $\displaystyle\f{H(M+h(p),p)-H(h(p),p)}{M} \leq -B$. Multiplying this inequality by $M^2 \geqslant 0$ we get $M H(M+h(p)) \leq -  B M^2$ because $H(h(p),p)=0$ for all $p \in [0,1]$.

Now, we bound each one of these terms of the right hand side of \eqref{eqM2} separately (keeping in mind the fact that $\epsilon$ will be chosen small enough)
\begin{equation}
\label{M1}
- B \int_{\R^d} M^2 \big( \f{1}{\epsilon} - M - h(p) \big) dx \leq - B \big( \f{1}{\epsilon} - h_0 - \tilde{K}' \big) \int_{\R^d} M^2 dx,
\end{equation}
where 
\begin{equation}
\label{def:K'}
\tilde{K}' = \sup \{ |n-h(p)|, \mbox{ as } |n|\leq \sup_{0<\epsilon\leqslant\epsilon_0} \|n\|_{L^\infty} \text{ and } 0 \leqslant p \leqslant 1\}.
\end{equation} 
From Lemma \ref{infbounds}, $\tilde{K}'$ is finite and by definition 
$M=n-h(p)$, $\tilde{K}'$ bounds $|M|$.
We pick $\epsilon < \epsilon_0$ such that $1 - \epsilon n > \frac{1}{2}$ (again, using Lemma~\ref{infbounds}). After using a Cauchy-Schwarz inequality, we get
\begin{multline}
\label{M2}
\left\lvert \int_{\R^d} 2 M h'(p) \epsilon\f{1}{1 - \epsilon n} (A(x) \nabla M + h'(p) A(x) \nabla p) \cdot \nabla p  dx \right\rvert\\
\leqslant h'_0 \tilde{K}' \epsilon \Big( (h'_0 + \frac{1}{2}) \int_{\R^d} \nabla p \cdot (A(x) \nabla p) dx  + \frac{1}{2} \int_{\R^d} \nabla M \cdot (A(x) \nabla M) dx \Big).
\end{multline}
Finally, by definition of $H$ in \eqref{def:H}, we have
$$
(-F_2-H)(n,p) = p(F_1-F_2)(n,p).
$$
Then using the constant $K_F$ defined in \eqref{def:KF}, we deduce, applying 
a Cauchy-Schwarz inequality,
\begin{multline}
\label{M3}
\left\lvert \int_{\R^d} \big( M h''(p) \nabla p \cdot ( A(x) \nabla p ) - M h'(p) (1-p) (- F_2 - H) (M + h(p), p) \big) dx \right\rvert \\
\leqslant h''_0 \tilde{K}' \int_{\R^d} \nabla p \cdot (A(x) \nabla p) dx + \frac{h'_0 K_{F}}{2} \int_{\R^d} M^2 dx + \frac{h'_0K_{F}}{2} \int_{\R^d} p^2 dx.
\end{multline}
Combining \eqref{M1}, \eqref{M2} and \eqref{M3} we get
\begin{multline*}
\displaystyle\frac{1}{2} \displaystyle\frac{d}{dt} \int_{\R^d} M^2 dx + \int_{\R^d} \nabla M \cdot (A(x) \nabla M) dx 
\\ \leqslant - B (\f{1}{\epsilon} - h_0 - \tilde{K}' - \frac{h'_0 K_F}{2 B}) \int_{\R^d} M^2 dx  + \frac{h'_0 K_F}{2} \int_{\R^d} p^2 dx 
\\ + \big(\tilde{K}' h''_0 + h'_0 \epsilon (h'_0 + \frac{1}{2}) \big) \int_{\R^d} \nabla p \cdot (A(x) \nabla p) dx 
\\ + h'_0 \frac{\epsilon}{2} \int_{\R^d} \nabla M \cdot (A(x) \nabla M )dx.
\end{multline*}
This is the expected estimate.
\end{proof}

With Lemmas \ref{lem:pestimate} and \ref{lem:Mestimate} we can proceed to prove Proposition \ref{lem:universalEstimates}.
\begin{proof}[Proof of Proposition \ref{lem:universalEstimates}]
Let $\alpha > 0$, summing the inequality in Lemmas \ref{lem:pestimate} and \ref{lem:Mestimate},
we obtain
\begin{multline*}
\f{1}{2} \f{d}{dt} \int_{\R^d} (M^2 + \alpha p^2) dx + (1 - \epsilon C_3 - \alpha \epsilon C_2) \int_{\R^d} \nabla M \cdot (A(x) \nabla M) dx \\+ ( \alpha(1 - \epsilon C_1) - C_6) \int_{\R^d} \nabla p \cdot (A(x) \nabla p) dx 
\leq \big(C_4 - \f{C_5}{\epsilon} \big) \int_{\R^d} M^2 dx + \big(C_6 + \alpha K \big) \int_{\R^d} p^2 dx.
\end{multline*}
Now, we can pick $\alpha > 0$, $\epsilon_0'\in (0,\epsilon_0)$ such that for all $0 < \epsilon < \epsilon_0'$,
$$
1 - \epsilon (C_3 + \alpha C_2) \geq \f{1}{2}, \quad \mbox{ and } \quad
\alpha (1 - \epsilon C_1) - C_6 \geq \f{1}{2}.
$$
The choice $\alpha = C_6 + 1$, then $\epsilon_0' = \min(\f{1}{2 (C_3 + C_2 (C_6 + 1))} ,\f{1}{2 C_1 (C_6 + 1)})$ suffices.
Hence we arrive at
\begin{multline}\label{eqMp}
\f{d}{dt} \int_{\R^d} (M^2 + \alpha p^2) dx + \int_{\R^d} \nabla M \cdot (A(x) \nabla M) dx + \int_{\R^d} \nabla p \cdot (A(x) \nabla p) dx  \\
\leq 2\big(C_4 - \f{C_5}{\epsilon} \big) \int_{\R^d} M^2 dx + 2\big(C_6 + \alpha K \big) \int_{\R^d} p^2 dx.
\end{multline}

Next, using the positivity of $A$, we may write for all $\epsilon > 0$ smaller than $\epsilon_0$ and $C_4 / C_5$
\[
	\f{d}{dt} \int_{\R^d} (M^2 + \alpha p^2) dx \leq 2 \f{C_6 + \alpha K}{\alpha} \int_{\R^d} (\alpha p^2 + M^2) dx,
\]
and thus by Gronwall's lemma, for all $\epsilon > 0$ small enough, with $C_0 := 2 \f{C_6 + \alpha K}{\alpha}$
\[
	\int_{\R^d} \big( M^2 (t, x) + \alpha p^2 (t, x) \big) dx \leq e^{C_0 t} \big( \lVert M^{\text{init}} \rVert_{L^2 (\R^d)}^2 + \alpha \lVert p^{\text{init}} \rVert_{L^2 (\R^d)}^2 \big).
\]
Since initial data are uniformly bounded in $L^2(\R^d)$ thanks to \eqref{hyp:init} and Remark \ref{rmkinit},
the first part of Proposition~\ref{lem:universalEstimates} is proved.
For all $T > 0$, $M$ and $p$ are uniformly bounded in $\mathcal{H}^1_T$ for $\epsilon$ small enough.

The second part follows easily from a time integration of \eqref{eqMp}. If $\epsilon$ is small enough, we get
\begin{multline*}
\int_0^t \int_{\R^d} \big( \nabla M \cdot( A(x) \nabla M)  + \nabla p \cdot( A(x) \nabla p) \big) dx ds \leq 2 \big(C_6 + \alpha K \big) \int_0^t \int_{\R^d} p^2 dx ds \\+ \int_{\R^d} \big((M^{\text{init}})^2 + \alpha (p^{\text{init}})^2 \big) dx.
\end{multline*}
Since we have proved the uniform $L^2$-bound of $p$, we conclude from the positivity of $A$ 
(Assumption \ref{assa}), and the uniform bounds on the initial data.
\end{proof}

\subsection{Convergence of $M$}

Until now we have not used the strength of the negative term in $\f{1}{\epsilon}$ in the right hand side of \eqref{eqMp}. Thanks to it, we can even get convergence of $M$.
\begin{lemma}
\label{lem:Mto0}
Under the assumptions of Theorem \ref{bigthm}, for all $T > 0$, 
$M \xrightarrow[\epsilon \to 0]{} 0$ strongly in $\mathcal{H}^1_T$.
\end{lemma}
\begin{proof}
Back to the estimate in Lemma \ref{lem:Mestimate}, and thanks to Proposition \ref{lem:universalEstimates}, we may write
\[
	\f{d}{dt} \int_{\R^d} M^2 dx \leq \big(2 C_4 - \f{2 C_5}{\epsilon}\big) \int_{\R^d} M^2 dx + C (t),
\]
where $C(t) := 2 C_6 \big( \int_{\R^d} p^2 dx + \int_{\R^d} \nabla p \cdot ( A(x) \nabla p ) dx)$.
From Proposition \ref{lem:universalEstimates}, we deduce that $C$ is bounded in $L^1(0,T)$.
Applying a Gronwall's lemma, we may write
\[
	\int_{\R^d} M^2 dx \leq e^{- 2 (\f{C_5}{\epsilon} - C_4) t}\big(\lVert M^{\text{init}} \rVert_{L^2 (\R^d)} + \int_0^t e^{2 (\f{C_5}{\epsilon} - C_4) t'} C(t') dt' \big).
\]
Let $\epsilon$ be small enough such that $\f{C_5}{\epsilon} > C_4$. 
Then integrating the latter inequality for $t\in[0,T]$, we deduce
$$
\int_0^T \int_{\R^d} M^2 dxdt \leq \frac{\epsilon}{2(C_5-\epsilon C_4)} \|M^{\text{init}}\|_{L^2(\R^d)} + \int_0^T \int_0^t e^{2 (\f{C_5}{\epsilon} - C_4) (t'-t)} C(t') dt'.
$$
We make a change of variable to estimate the last term in the right hand side:
$$
\begin{array}{ll}
\displaystyle \int_0^T \int_0^t e^{2 (\f{C_5}{\epsilon} - C_4) (t'-t)} C(t') dt'dt & 
\displaystyle = \int_0^T \int_{t'}^T e^{2 (\f{C_5}{\epsilon} - C_4) (t'-t)} dt \, C(t') dt'  \\[2mm]
& \displaystyle = \int_0^T \int_{t'-T}^0 e^{2 (\f{C_5}{\epsilon} - C_4) \tau} d\tau \, C(t') dt' \\[2mm]
& \displaystyle \leq \frac{\epsilon}{2(C_5-\epsilon C_4)} \int_0^T C(t') dt'.
\end{array}
$$
We conclude that 
$$
\int_0^T \int_{\R^d} M^2 dxdt \leq \frac{\epsilon}{2(C_5-\epsilon C_4)}
\left(\|M^{\text{init}}\|_{L^2(\R^d)} + \int_0^T C(t') dt'\right).
$$
It implies the expected convergence as $\epsilon\to 0$.
\end{proof}

\subsection{Compactness result and proof of Theorem \ref{bigthm}}

Before proving our main result, we recall the following compactness result (see \cite{Sim.Compact}).
\begin{nnlemma}[Lions-Aubin]
Let $T > 0, q \in (1, \infty)$, $(\psi_n)_n$ a bounded sequence in $L^q (0, T ; H)$, where $H$ is a Banach space. 
If $\psi_n$ is bounded in $L^q (0, T ; V)$ and $V$ compactly embeds in $H$, and if $(\partial_t \psi_n)_n$ is bounded in $L^q (0, T ; V')$ uniformly with respect to $n$, 
then $(\psi_n)_n$ is relatively compact in $L^q (0, T ; H)$.
\end{nnlemma}

\begin{proof}[Proof of Theorem \ref{bigthm}]
We split the proof into three steps. First, our previous estimates together with Lions-Aubin lemma enable us to prove relative compactness on bounded domains. 
Then, through a diagonal extraction process, we prove that there exists (up to extracting a subsequence) a global limit. 
Finally, thanks to our uniform estimates, we prove that this limit satisfies a universal equation whose solution is unique, which in turn implies convergence of the whole sequence.

\paragraph{Step 1: Local relative compactness.}
 For $R>0$ we define the increasing sequence $(B_R)_R$ of balls of radius $R$ with center $0$ in $\R^d$, and
 $H_R = L^2 (B_R)$, $V_R = H^1 (B_R) \cap L^{\infty} (B_R)$, and pick $T > 0$.
 Then, we check that Lions-Aubin Lemma with $q=2$ can be applied to
 $$
 \psi^{(R)}_{\epsilon} = \restriction{p^{\epsilon}}{B_R}.
 $$
 Lemma \ref{lem:pestimate} gives boundedness in $L^q(0, T ; V_R)$. The compact embedding is classical (Rellich-Kondrachov).
 We check that the time derivative is bounded. Let $\chi \in V_R$, $\langle \cdot, \cdot \rangle = \langle \cdot, \cdot \rangle_{V'_R, V_R}$ and $t \in (0, T)$.
 \begin{multline*}
  \int_0^t \lvert \langle \p_t p^{\epsilon} (\tau) , \chi \rangle \rvert^2 d\tau
  = \int_0^t \big\lvert \langle  \nabla \cdot (A(x) \nabla p^{\epsilon}) - 2 \epsilon \nabla p^{\epsilon} \cdot \f{A(x) \nabla n^{\epsilon}}{1 - \epsilon n^{\epsilon}} 
  + p^{\epsilon} (1 - p^{\epsilon}) (F_1 - F_2) (n^{\epsilon}, p^{\epsilon} ) , \chi \rangle \big\rvert^2 d\tau.
 \end{multline*}
This can be bounded
\begin{align*}
  \int_0^t \lvert \langle \p_t p^{\epsilon} (\tau) , \chi \rangle \rvert^2 d\tau &\leq \Big( \int_0^t \int_{B_R} \lvert A(x) \nabla p^{\epsilon} \cdot \nabla \chi \rvert  \Big)^2 \\
  &+\epsilon \Big(\int_0^t \int_{B_R} \lvert A(x) \nabla n^{\epsilon} \cdot \nabla p^{\epsilon} \rvert \Big)^2 + C t \int_{B_R} \chi^2 \\
  &\leq \lVert \nabla \chi \rVert^2_{H_R} \lVert \nabla p^{\epsilon} \rVert^2_{L^2(0, T ; H_R)}\\
  & + 2 \epsilon \lVert \nabla n^{\epsilon} \rVert^2_{L^2 (0, T ; H_R)} \lVert \nabla p^{\epsilon} \rVert^2_{L^2(0, T ; H_R)} \lVert \chi \rVert_{\infty}^2 + C T \lVert \chi \rVert^2_{H_R},
\end{align*}
which gives the required bound, uniform in $0 < \epsilon < \epsilon_0$, for $\epsilon_0$ small enough. 
This holds thanks to Lemmas~\ref{lem:pestimate}~and~\ref{lem:Mestimate}.

\paragraph{Step 2: Global convergence.} 
Now, for all $R \in \N$, one can extract converging (in $L^2 (0, T; H_R)$) 
subsequence from $(p^{\epsilon})_{\epsilon}$ by Lions-Aubin Lemma.
We perform a diagonal extraction process successively in $R$, so that
$$
p^{\epsilon_m^{(R)}} \xrightarrow[m \to \infty]{} p^{(R)} \text{ in } L^2 (0, T ; H_R),
$$
and by construction $(\epsilon_m^{(R_1)})_m$ is a subsequence of $(\epsilon_m^{(R_2)})_m$ if $R_2 > R_1$.
Because the whole family $(p^{\epsilon})_{\epsilon}$ is in $L^2(0, T ; H^1 (\R^d))$ 
uniformly in $\epsilon$ (by Lemma \ref{lem:pestimate}),
one gets weak convergence of gradient
$$
\nabla p^{\epsilon_m^{(R)}} \xrightharpoonup[m \to \infty]{} \nabla p^{(R)} \text{ in } L^2 (0, T ; H_R).
$$
Thanks to Lemma \ref{lem:pestimate}, we know that the limits $p^{(R)}$ are well-defined, do not depend on the extracted subsequences, satisfy the same bounds as $(p^{\epsilon})_{\epsilon}$ and
$$
R_2 > R_1 \implies \restriction{p^{(R_2)}}{B_{R_1}} = p^{(R_1)}.
$$
Therefore we can define $p^{0} \in L^2(0, T ; L^2 (\R^d))$ and we have constructed
a subsequence, still denoted $(p^\epsilon)_\epsilon$, such that 
$p^\epsilon \xrightarrow[\epsilon \to 0]{} p^0$ strongly in $L^2(0,T;L^2(B_R))$ for all
$R>0$.

To pass from local to global convergence, we need to have uniform in $\epsilon$ 
estimate in the tails $|x|>R$. To do so, let us introduce $\phi\in \mathcal{C}^\infty(\R^d)$
such that $0\leq \phi\leq 1$, $\phi(x)=0$ if $|x|<1/2$ and $\phi(x)=1$ if $|x|>1$.
Then we denote $\phi_R(x)=\phi(x/R)$.
Multiplying the equation satisfied by $p^\epsilon$ in \eqref{N1p} by $p^\epsilon \phi_R$ 
and integrating over $\R^d$, we deduce
$$
\begin{array}{r}
\displaystyle \frac 12 \frac{d}{dt} \int_{\R^d} (p^\epsilon)^2 \phi_R \,dx + \int_{\R^d} \nabla(p^\epsilon \phi_R)\cdot A(x)\nabla p^\epsilon \,dx + \int_{\R^d} \frac{2\epsilon}{1-\epsilon n^\epsilon} \phi_R p^\epsilon \nabla p^\epsilon \cdot A(x) \nabla n^\epsilon \,dx \\[2mm]
\displaystyle \leq K_F \int_{\R^d} (p^\epsilon)^2 \phi_R \,dx,
\end{array}
$$
where $K_F$ has been defined in \eqref{def:KF}. Using a Cauchy-Schwarz inequality,
we have
$$
\begin{array}{rcl}
\displaystyle \int_{\R^d} \nabla(\phi_R p^\epsilon)\cdot A(x)\nabla p^\epsilon \,dx & = &
\displaystyle  \int_{\R^d} p^\epsilon \nabla\phi_R \cdot A(x)\nabla p^\epsilon \,dx +
\int_{\R^d} \phi_R \nabla p^\epsilon\cdot A(x)\nabla p^\epsilon \,dx \\[2mm]
&\geq & \displaystyle 
- \left(\int_{\R^d} \nabla \phi_R \cdot A(x) \nabla \phi_R \,dx \right)^{1/2}
\left(\int_{\R^d} \nabla p^\epsilon \cdot A(x) \nabla p^\epsilon \,dx \right)^{1/2}.
\end{array}
$$
By definition of $\phi_R$ we have that $\nabla \phi_R(x) = \frac{1}{R} \nabla \phi(x/R)$
and $\nabla \phi_R(x)=0$ on $B_{R/2}\cup \R^d\setminus \overline{B_R}$.
As above, we take $\epsilon$ small enough such that $1-\epsilon n^\epsilon\geq \frac 12$,
which can be done thanks to Lemma \ref{infbounds}. 
Then, as in the proof of Proposition \ref{lem:universalEstimates}, there exists 
a nonnegative function 
$C(t)\in L^1(0,T)$ such that, thanks to a Cauchy-Schwarz inequality 
$$
\left|\int_{\R^d} \frac{2\epsilon}{1-\epsilon n^\epsilon} \phi_R p^\epsilon \nabla p^\epsilon \cdot A(x) \nabla n^\epsilon \,dx \right | \leq C(t) \epsilon,
\qquad \mbox{ and }\quad  
-\int_{\R^d} \nabla \phi_R\cdot A(x)\nabla p^\epsilon \,dx \leq \frac{C(t)}{R}.
$$
Then, we have obtained
$$
\frac 12 \frac{d}{dt} \int_{\R^d} (p^\epsilon)^2 \phi_R \,dx \leq C(t)\big(\frac{1}{R} + \epsilon\big) +
K_F \int_{\R^d} (p^\epsilon)^2 \phi_R \,dx.
$$
Using a Gronwall Lemma, it implies
$$
\int_{\R^d} (p^\epsilon)^2 \phi_R \,dx \leq  e^{2 K_F t} \int_{\R^d} (p^{\text{init}})^2 \phi_R \,dx+ 
\big(\frac{1}{R} + \epsilon\big) \int_0^t 2C(\tau) e^{2K_F(t-\tau)}\,d\tau.
$$
By definition of $\phi_R$ ($\phi_R(x)=1$ on $\R^d\setminus B_R$),
we deduce that for all $\epsilon>0$ small enough and all $R>0$,
\begin{equation}\label{eq:tail}
\int_0^T\int_{\R^d\setminus B_R} |p^\epsilon|^2 \,dx \leq \int_0^T\int_{\R^d} (p^\epsilon)^2 \phi_R \,dx 
\leq \frac{e^{2K_F T}-1}{2K_F} \left( \int_{\R^d\setminus B_{R/2}} (p^{\text{init}})^2 \,dx  
+ (\frac{1}{R} + \epsilon) \int_0^T 2C(t)\,dt\right).
\end{equation}
It implies an uniform bound, since $p^{\text{init}}\in L^2(\R^d)$.

Finally, we conclude that the subsequence $(p^\epsilon)_\epsilon$ converges strongly
towards $p^0$ in $L^2(0,T;L^2(\R^d))$ as $\epsilon\to 0$. Indeed, we have
$$
\int_0^T \int_{\R^d} |p^\epsilon - p^0|^2 \,dxdt = \int_0^T \int_{B_R} |p^\epsilon - p^0|^2 \,dxdt
+\int_0^T \int_{\R^d\setminus B_R} |p^\epsilon - p^0|^2 \,dxdt.
$$
The second term of the right hand side is uniformly bounded for $R$ large enough thanks to
\eqref{eq:tail} and the fact that $p^0\in L^2(0,T;L^2(\R^d))$.
For the first term we use the local convergence.

\paragraph{Step 3: Limit equation.} 
From the strong convergence of the sequence $(p^{\epsilon})\epsilon$ in $L^2(0,T;L^2(\R^d))$
and the Lipschitz continuity of the function $h$, we deduce that $(h(p^\epsilon))_\epsilon$
converges strongly in $L^2(0,T;L^2(\R^d))$ towards $h(p^0)$.
Moreover, using the triangle inequality, we have
$$
|n^\epsilon - h(p^{0}) | \leq |n^\epsilon - h(p^{\epsilon})| + h'_0 |p^\epsilon-p^0|.
$$
Applying Lemma \ref{lem:Mto0}, we deduce that 
\begin{equation}
 \label{limiteqN1}
n^\epsilon \xrightarrow[\epsilon \to 0]{} n^{0} := h (p^{0})
\text{ strongly in } L^2(0,T;L^2(\R^d)).
\end{equation}

Then, we obtain the equation satisfied by $p^0$ 
using the weak forms of the equations on $p^{\epsilon}$ in \eqref{N1p}: 
for all $\chi \in \mathcal{C}_c^{\infty} (\R^d)$,
\begin{align*}
\int_{\R^d} p^{\epsilon} (T, x) \chi(x) dx 
- \underbrace{\int_{\R^d} p^{\text{init},\epsilon} (x) \chi(x) dx}_{\text{weak convergence}} 
+ \underbrace{\int_0^T \int_{\R^d} \nabla p^{\epsilon}(t,x) \cdot A (x) \nabla \chi(x) dx dt}_{\text{weak convergence}}
\\ + 2 \epsilon \underbrace{\int_0^T \int_{\R^d} \chi(x) \nabla p^{\epsilon}(t, x) \cdot \f{A(x) \nabla n^{\epsilon}(t,x)}{1 - \epsilon n^{\epsilon}(t, x)} dx dt}_{\text{bounded as } \epsilon \to 0}
\\ = \int_0^T \int_{\R^d} \chi(x) \underbrace{p^{\epsilon} (1 - p^{\epsilon}) (F_1 - F_2)(n^{\epsilon}, p^{\epsilon})}_{\text{strong convergence}} dx dt
\end{align*}
We can pass to the limit in each term, using also \eqref{hyp:init} for the second term.

Hence $p^{0}$ is in $L^2(0,T;H^1(\R^d))$ 
and is a weak solution of the initial value problem
\begin{equation}
\label{limitsysp}
 \begin{cases}
    &\p_t p^0 - \nabla \cdot (A(x) \nabla p^0 = p^0 (1 - p^0) (F_1 - F_2)(n^0, p^0), \\
    & p^0(t=0, \cdot) = p^{\text{init}}.
 \end{cases}
\end{equation}
Using \eqref{limiteqN1} in \eqref{limitsysp} yields a self-contained initial valued 
reaction-diffusion system on $p^0$ that has a unique solution.
It defines in turn uniquely $n^0$ through \eqref{limiteqN1}.
Since solutions to the initial value system \eqref{limitsysp} are unique, 
all extracted subsequences converge to the same limit. 
Therefore, the whole sequences converge, strongly in $L^2$ with weak convergence of gradients.

This concludes the proof of Theorem \ref{bigthm}.

\end{proof}

\section{Generalization of the result}
\label{generalization}
We have stated Theorem \ref{bigthm} so as to keep simplicity and stick to the biological application in Section \ref{exwolbachia}. 
It can be slightly generalized 
in order to encompass spontaneous transition between variants.

Individuals in state $1$ may give birth to individuals in state $2$, and \textit{vice versa}.
To do so, we consider more general reaction term and replace system \eqref{syst1} by
$$
\begin{cases}
  \partial_t n_1 - \nabla \cdot (A(x) \nabla n_1) &= \widetilde{f}_1 (n_1, n_2),\\
  \partial_t n_2 - \nabla \cdot (A(x) \nabla n_2) &= \widetilde{f}_2 (n_1, n_2),
\end{cases}
$$

In fact, the basic property we require in our proof is that $p$ stays between $0$ and $1$, that is, $n_i$ remain non-negative.
Here is the minimal hypothesis ensuring positivity (in the spirit of \cite{Perthame}).
\begin{assumption}[Positivity]
\label{asspos}
We assume
$$ \forall n_1, n_2 \in \R_+, \quad \widetilde{f}_1(0, n_2) \geq 0 \mbox{ and } \widetilde{f}_2(n_1, 0) \geq 0.$$
\end{assumption}
\begin{proof}[Proof of ``Assumption \ref{asspos} implies positivity``]
We prove that if the initial data $n_1^{\text{init}}$, $n_2^{\text{init}}$ are non-negative and if Assumption \ref{asspos} holds, then $n_1$ and $n_2$ remain non-negative.
It is a simple application of the comparison principle for this parabolic system.
A solution that lies initially above a sub-solution remains above it. The constant $(0, 0)$ is indeed a sub-solution.
\end{proof}

For the sake of clarity of the presentation, we only consider an extension 
of the biological example from Section \ref{exwolbachia}.
This allows us to take into account imperfect maternal transmission.
We assume that at a rate $\mu$, infected females lay eggs which do not carry \textit{Wolbachia}.
This quantity is very commonly tested by entomologists, and usually shown to be close to $0$ (see \cite{Wer.Wolbachia} and references, and for example \cite{Dut.From} where they obtained $\mu = 0.04$ and $\mu = 0$).
This feature is included in the following model taken from \cite{Fen.Solving} 
(neglecting the pathogen effect),
\begin{equation}
 \label{eq:Fen}
 \begin{cases}
    &\p_t n_i - \nabla \cdot (A(x) \nabla n_i) = n_i F_u (1-s_f) (1 - \mu) - n_i(d_i + \sigma(n_i + n_u)),\\
    &\p_t n_u - \nabla \cdot (A(x) \nabla n_u) = n_u F_u (1 - s_h\frac{n_i}{n_u + n_i}) + \mu F_u (1-s_f) n_i - n_u (d_u + \sigma(n_i + n_u)).
 \end{cases}
\end{equation}
Here, the reduced population would be $n = \sigma (n_i + n_u)$.
The corresponding dynamics in $(n, p)$ reads,
\begin{equation}
\label{eq:FenNp}
 \begin{cases}
    &\p_t n - \nabla \cdot (A(x) \nabla n) = n \Big( F_u \big( p (1-s_f) + (1-p) (1-s_h p) \big) - d_u \big( (\delta - 1)p +1\big) - n \Big),\\
    &\p_t p - \nabla \cdot (A(x) \nabla p) - 2 \frac{\nabla n}{n} A(x) \nabla p = p \Big( (1-p) \big(F_u (1-s_hp) - d_u (\delta - 1) \big) - \mu F_u (1-s_f) \Big).
 \end{cases}
\end{equation}
We notice in particular that the reaction term for $p$ in \eqref{eq:FenNp} does not depend on $n$.
It yields directly the equation \eqref{eq:BarTur} with a function $n$ in the left 
hand side that depends on $p$, whereas in \cite{BarTur.Spatial} the function $n$ in the
gradient in the left hand side is assumed to be given.

As in Section \ref{exwolbachia}, we introduce the parameter $\epsilon$ to characterize 
the high fertility and strong competition and propose the following extension of 
system \eqref{wolsys1}, with imperfect maternal transmission,
\begin{equation}\label{wolsyst2}
 \begin{cases}
    \p_t n_i - \nabla \cdot (A(x) \nabla n_i) = (1-\mu) (1-s_f) F_u n_i \big(\f{1}{\epsilon} - \sigma (n_i + n_u) \big)_+ - \delta d_u n_i, \\
    \p_t n_u - \nabla \cdot (A(x) \nabla n_u) = F_u \big(n_u (1 - s_h p) + \mu (1-s_f) n_i p \big) \big(\f{1}{\epsilon} - \sigma (n_i + n_u) \big)_+ - d_u n_u, \\
 \end{cases}
\end{equation}
with $p=\frac{n_i}{n_i+n_u}$ as usual.
In this system, the notation $a_+=\max\{0,a\}$ denotes 
the positive part of $a \in \R$.

For the reduction, as above, we identify $n_1=n_i$ and $n_2=n_u$ and we deduce 
from \eqref{wolsyst2} the equations satisfied by $n=\f{1}{\epsilon}-\sigma(n_i + n_u)$ 
and $p$,
\begin{equation}\label{wolsyst2n}
\begin{array}{rl}
\displaystyle \p_t n - \nabla \cdot (A(x)\nabla n) = & \displaystyle
- \big( \f{1}{\epsilon} -n \big) F_u\left((1-s_f) ((1-\mu) p +\mu p^2)+(1-p)(1-s_h p)\right)n_+
\\[2mm]
& \displaystyle + d_u (\delta p + 1-p) \big( \f{1}{\epsilon} -  n \big) ,\\[2mm]
\end{array}
\end{equation}
\begin{equation}\label{wolsyst2p}
\begin{array}{rl}
\displaystyle  \p_t p - \nabla \cdot (A(x) \nabla p) +2 \frac{\nabla n}{n} A(x) \nabla p = 
& \displaystyle F_u p \Big( (1-p)\big((1-\mu)(1-s_f) - (1-s_h) p\big)   \\[2mm]
& \displaystyle + \mu (1-s_f) p^2 \Big) n_+ + p(1-p)d_u (1-\delta).
\end{array}
\end{equation}
Using the notation in \eqref{def:Nnp}, we define as in \eqref{def:H} the function $H$ by
$$
\begin{array}{rcl}
 H (n, p) &:=& - F_u n \big( p (1- \mu) (1 - s_f) + (1 - p)(1-s_h p) + \mu (1-s_f) p^2\big) + d_u (p (\delta -1) + 1)  \\[2mm]
& = & - F_u  n \big((s_h + \mu(1-s_f) ) p^2 - (s_f + s_h + \mu(1-s_f)) p + 1\big) + d_u ((\delta -1)p + 1).
\end{array}
$$
When $\mu=0$, we notice that we recover the same expression as in the case of perfect 
maternal transmission in Section \ref{exwolbachia}.
Then, the function $h$ and the reaction term are modified. 

In this case, as in Lemma \ref{equilibria}, we may investigate the equilibria 
of \eqref{wolsyst2n}--\eqref{wolsyst2p}. We get from straightforward computations:
\begin{lemma}
Let 
\[
\Delta = \big(\delta(s_f + s_h) + (\delta - 1 - \mu)(1 - s_f) \big)^2 - 4 \delta \big(s_h + \mu(1-s_f) \big) \big( \delta - (1- \mu) (1-s_f) \big).
\]
Let us assume that $\Delta > 0$. 
When $\mu = 0$, the condition $\Delta > 0$ is equivalent to 
$\big(\delta s_h - \delta + (1-s_f) \big)^2 > 0$ which is always satisfied.
Then, there are 4 equilibria associated to the system 
\eqref{wolsyst2p}--\eqref{wolsyst2n} in the reduced variable $(n,p)$:
\begin{itemize}
\item The co-existence equilibrium reads
$$
\bepa
p^*_C = 1 - \displaystyle\f{\delta(s_f + s_h) + (\delta - 1 + \mu)(1-s_f) - \sqrt{\Delta}}{2 \delta (s_h+ \mu(1-s_f))},
\\[10pt]
n^*_C = \displaystyle\f{\delta d_u}{(1-\mu)(1-s_f) F_u}, 
\eepa
$$
it remains unstable.

\item The steady state $(0,0)$ is unstable.

\item The stable \textit{Wolbachia} invasion equilibrium reads
$$
\bepa
p^*_W = 1 - \displaystyle \f{\delta(s_f + s_h) + (\delta - 1 + \mu)(1-s_f) + \sqrt{\Delta}}{2 \delta (s_h+ \mu(1-s_f))} < 1,
\\[10pt]
n^*_W = \displaystyle\f{\delta d_u}{(1-\mu)(1-s_f) F_u} = n^*_C.
\eepa
$$

\item The stable \textit{Wolbachia} extinction equilibrium is unchanged:
$n_E^* = \displaystyle\f{d_u}{F_u}$, $p_E^* = 0.$
\end{itemize}

\end{lemma}

From straightforward computation, we may adapt Theorem \ref{bigthm}
in this framework. Then, the analogue of Corollary \ref{cor1} reads
\begin{corollary}
Assume that $A$ satisfies Assumption \ref{assa}.
Given $n_1^{\text{init},\epsilon}$ and 
$n_2^{\text{init},\epsilon}$ such that there exists $p^{\text{init}}\in L^2(\R^d)$ such that 
$p^{\text{init},\epsilon}\rightharpoonup p^{\text{init}}$ as $\epsilon\to 0$ in $L^2(\R^d)$-weak 
and $\f{1}{\epsilon} - \sigma(n_1^{\text{init}, \epsilon} + n_2^{\text{init}, \epsilon}) - \f{d_u}{F_u} \in L^2 \cap L^{\infty} (\R^d)$ with uniform bounds in $\epsilon > 0$,
then Theorem \ref{bigthm} applies and the solutions $(n_i^{\epsilon}, n_u^{\epsilon})_{\epsilon > 0}$ 
of \eqref{wolsyst2}
  satisfy the convergence result in \eqref{conv}. The limiting equation reads
 \begin{equation}
 \p_t p - \nabla \cdot (A(x) \nabla p) = r_{\mu}(p),
\label{limeqWolba2}
\end{equation}
where
\begin{align*}
r_{\mu} (p) &= d_u p \Big((1-\mu)(1-s_f) \f{(\delta - 1) p + 1}{(s_h + \mu(1-s_f)) p^2 - (s_f + s_h + \mu(1-s_f))p + 1}  - \delta \Big).
\end{align*}
\end{corollary}
For small $\mu$, $r_{\mu}$ is still a bistable function provided $\Delta > 0$,
however the stable state $1$ is displaced. 

\begin{figure}[h!]
 \includegraphics[width=.5\linewidth]{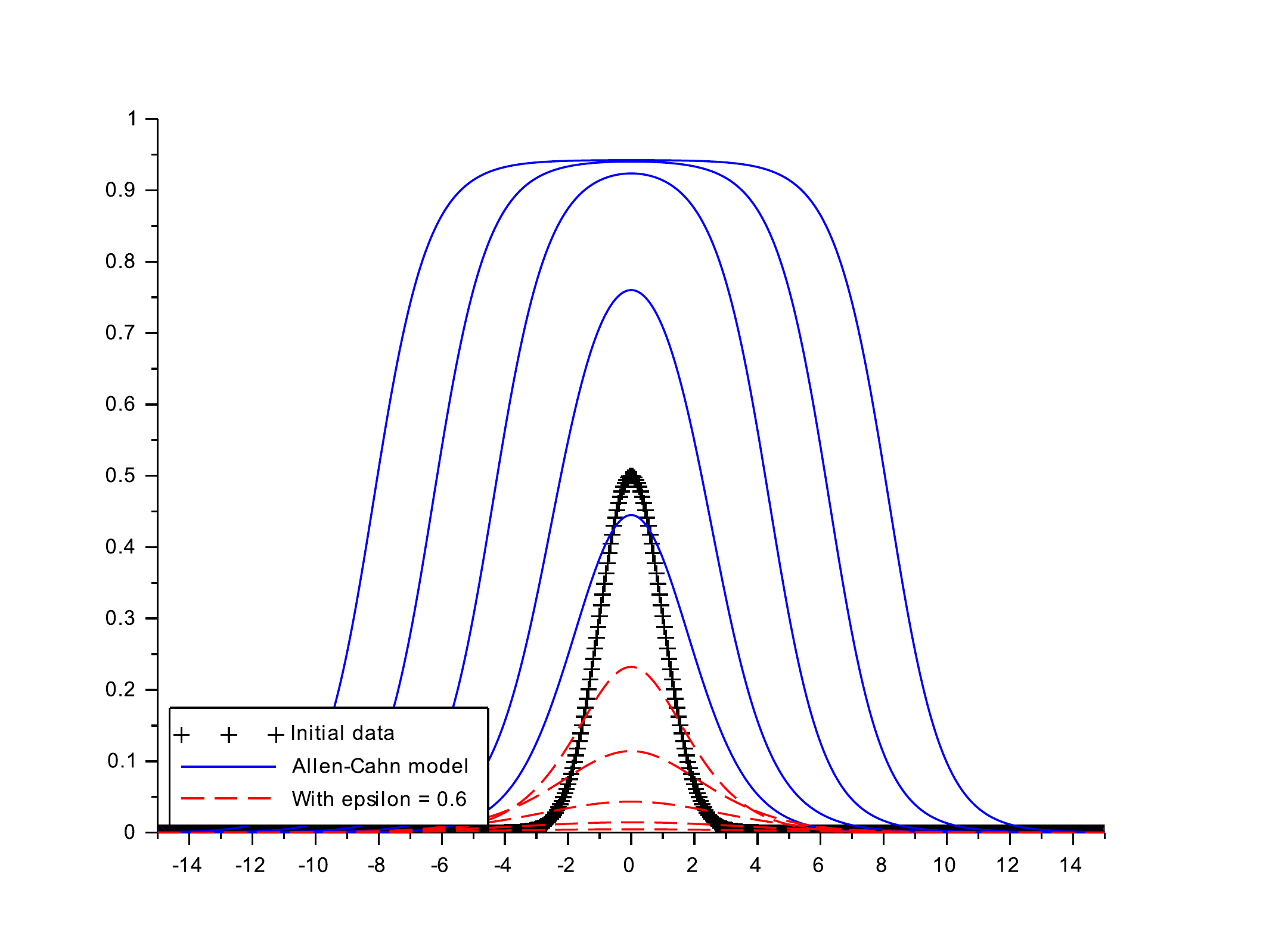}
 \includegraphics[width=.5\linewidth]{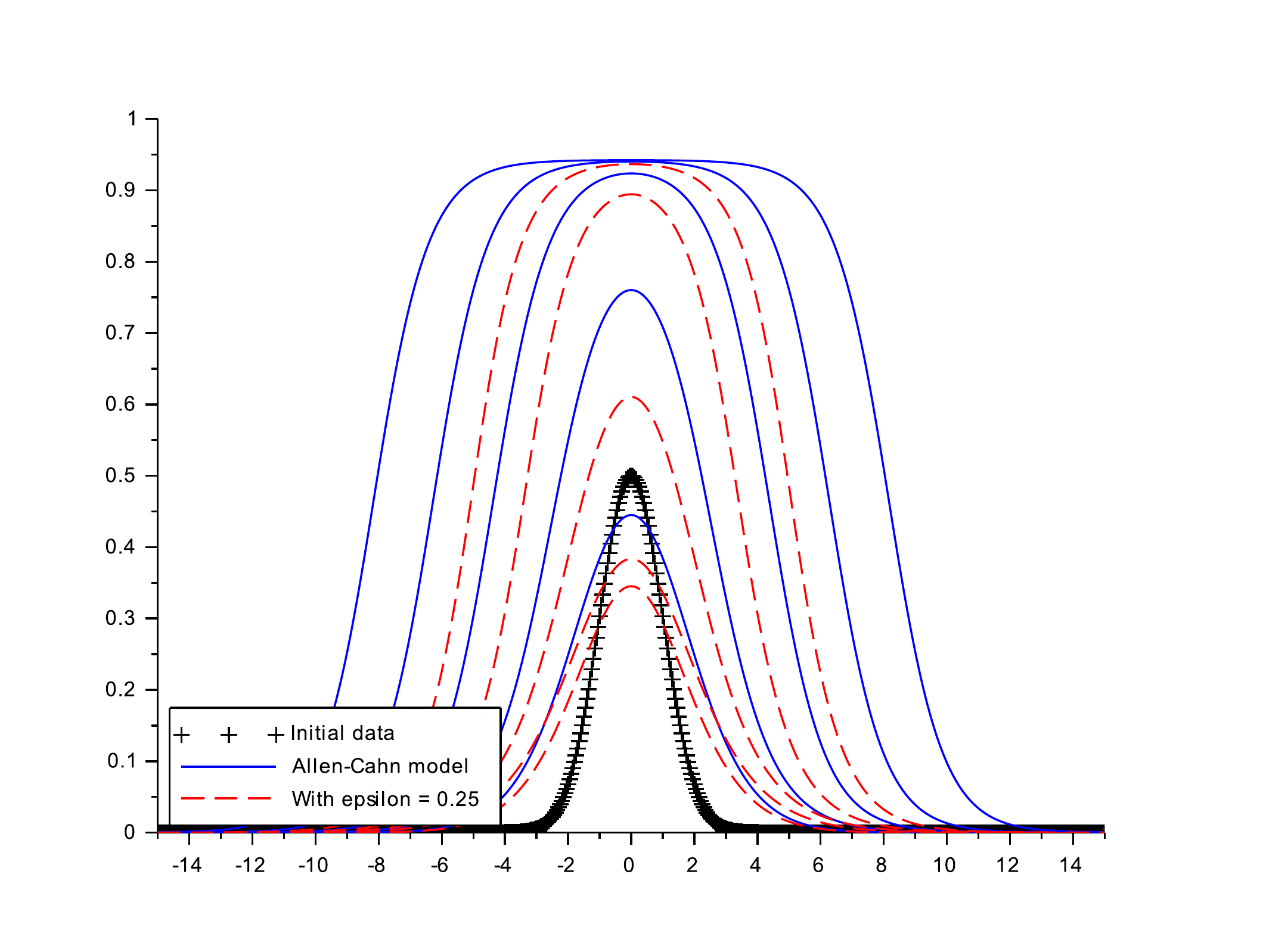}
 \includegraphics[width=.5\linewidth]{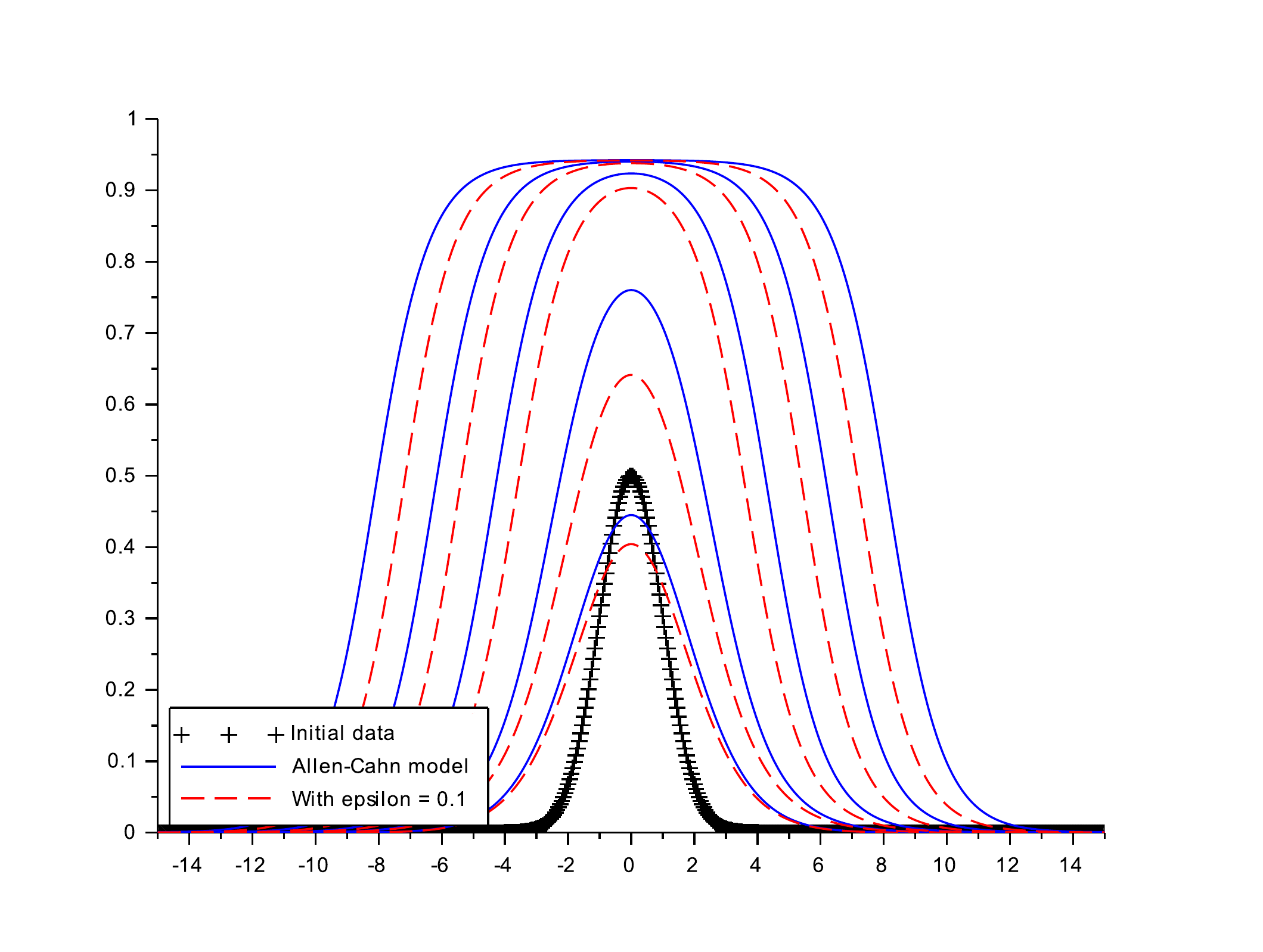}
  \includegraphics[width=.5\linewidth]{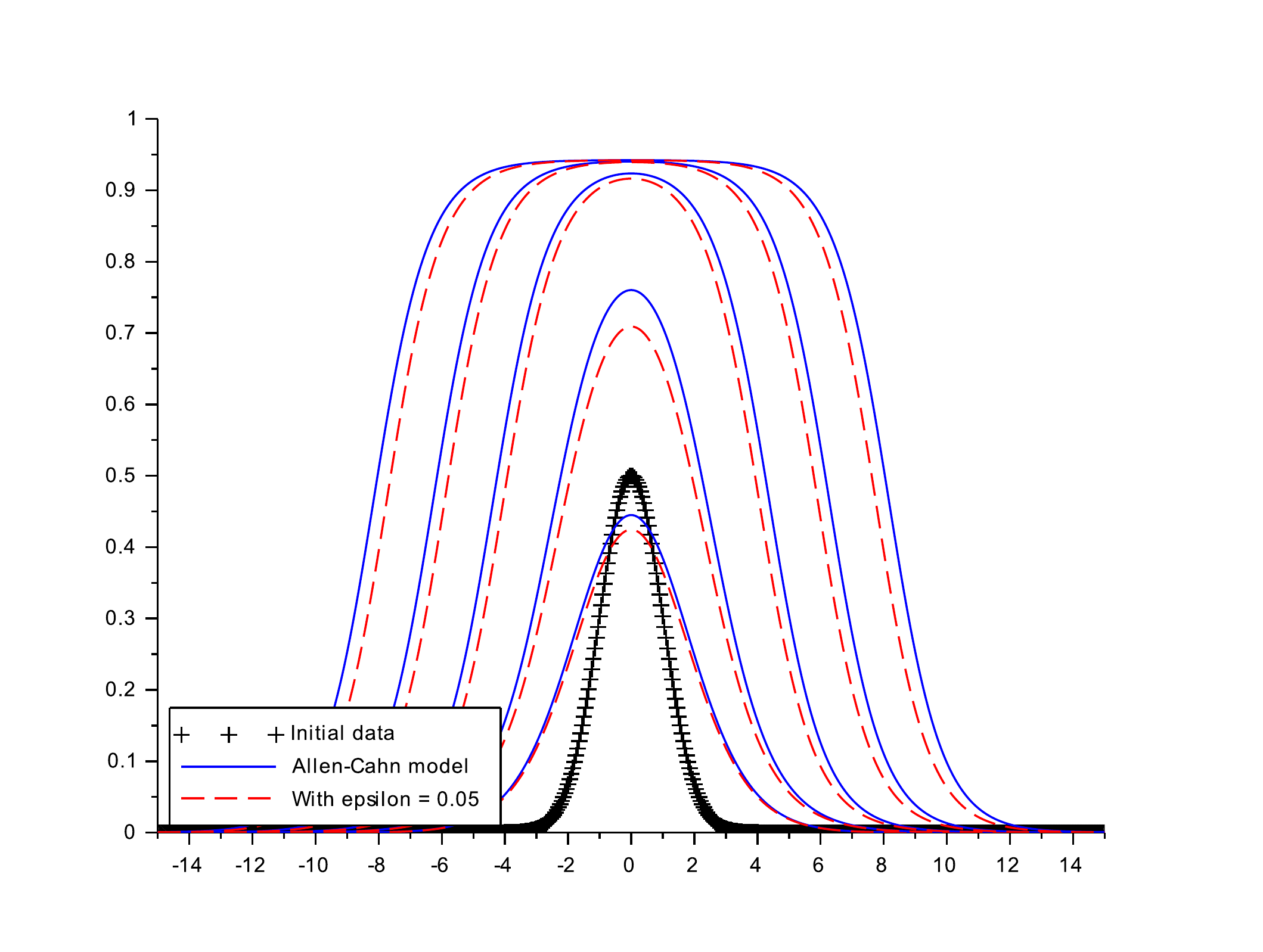}
 \caption{Initial data (+) creating a traveling wave in the limit system (blue) and convergence of the two-population solution (dashed red) as $\epsilon$ diminishes.}
 \label{fig:convergence2}
\end{figure}

We give a numerical illustration of this case in Figure \ref{fig:convergence2}, using
a similar approach as in Section \ref{sec:num}.
We use the same parameters as for Figure \ref{fig:convergence}, except that $s_f=0$, $\mu = .04$ and the initial data is smaller (less infected mosquitoes are introduced).

For Figure \ref{fig:convergence2}, we use the same discretization and numerical scheme as in Figure \ref{fig:convergence}.
The blue lines represent the solution of the limiting system \eqref{limeqWolba2}.
In dashed red lines are plotted the numerical results for the system of
two populations \eqref{wolsyst2p}.
We observe that the solution of the limiting bistable system \eqref{limeqWolba2}
exhibits a traveling front which propagates into the whole domain. 
Then the numerical results for 4 different values of the parameter $\epsilon$ 
are represented. 
For large populations, we observe that as $\epsilon$ goes to $0$ (recall that the order of magnitude of the population size is $\f{1}{\sigma \epsilon}$), 
the solution to the whole system \eqref{wolsyst2p} gets closer to the one of the limiting system.
However, for small populations, we see a clear modification of the wave's shape and speed, which is slower than the limit wave.

\section{Conclusion and perspectives}
\label{conclusion}

We have established in this paper the rigorous convergence, under suitable assumptions,
of a 2 by 2 reaction diffusion model of Lotka-Volterra type towards a simple model
for the frequency of a variant. It justifies the use of such reduced model 
in applications.
Let us discuss quickly our scaling choice in Assumption \ref{asseps}, 
in the case of \textit{Wolbachia}.

Another biologically relevant scaling assumption would not give a limiting system 
consisting in only one equation on frequency.
Indeed, if we consider the following alternative model
\begin{equation}
\label{wolsys2}
 \begin{cases}
    \p_t n_i - \nabla \cdot (A(x) \nabla n_i) &= (1-s_f) F_u n_i \big(1 - \epsilon \sigma (n_i + n_u) \big) - \delta d_u n_i, \\
    \p_t n_u - \nabla \cdot (A(x) \nabla n_u) &= F_u n_u (1 - s_h \f{n_i}{n_i + n_u} ) \big(1- \epsilon \sigma (n_i + n_u) \big) - d_u n_u. \\
 \end{cases}
\end{equation}
Then, $n$ and $p$ satisfy the following system, that does not depend on $\epsilon$
\begin{equation}
\label{wolsys2Np}
 \begin{cases}
    \p_t n - \nabla \cdot (A(x) \nabla n) = F_u \big(1 - n \big) \big( A(p) - B(p)  n \big), \\
    \p_t p - \nabla \cdot (A(x) \nabla p) + \f{2 \nabla p \cdot A(x) \nabla n}{1-n} = p (1-p) \big( F_u n (s_h p - s_f) - d_u (\delta - 1) \big), \\
 \end{cases}
\end{equation}
where
\begin{equation*}
 \begin{cases}
  A(p) &= ((\delta - 1)p + 1) \f{d_u}{F_u},\\
  B(p) &= s_h p^2 -(s_f + s_h) p + 1. \\
 \end{cases}
\end{equation*}
The dependancy in $\epsilon$ in the resulting model is only through the initial data.
Thus, $\epsilon \to 0$ does not imply $n - \f{A(p)}{B(p)} \to 0$ in \eqref{wolsys2}, \eqref{wolsys2Np}.


We conclude that the use of simple bistable models for the spatial spread 
of \textit{Wolbachia} can be justified mathematically.
This is the object of Theorem \ref{bigthm}.
However, we must keep in mind that this result applies 
only if population size and fecundity scale properly.

In the context of \textit{Wolbachia} modeling, bistable equations like \eqref{eq:BarTur} have been used (for example in \cite{BarTur.Spatial} or \cite{Sch.Constraints}) because they provide with a unique (up to translations) and linearly stable traveling wave solution.
Hence, with a bistable model at hand we can compute a speed that may be interpreted as an invasion speed.

Therefore a natural continuation of the present work would be to try and specify Theorem \ref{bigthm} to traveling waves. 
The open question reads: does the frequency in the two-populations
model converge to the unique traveling wave solution of the limit bistable equation?
If yes, in what sense?
Indeed, there are two types of convergence involved: on the first hand in the singular limit (where we identified
a small parameter $\epsilon$), that proves convergence of the system’s frequency to a solution of the limit
bistable equation; and on the other hand the well-known attractiveness result of the unique traveling wave solution
in the bistable case.
Moreover, existence and (local) stability of traveling waves has been proved for competitive systems (see \cite{Gardner} for example). 
How to compare the traveling speed for competitive system with the one for the reduced model on the frequency? 

\bigskip

{\bf Acknowledgements.}
The authors acknowledge partial supports from the Capes/Cofecub project
Ma-833 15 {\it ``Modeling innovative control method for Dengue fever''}
and from the Programme Convergence Sorbonne Universit\'es / FAPERJ {\it ``Control 
and identification for mathematical models of Dengue epidemics''}.

They warmly thank B. Perthame for his patient and constant help, useful discussions
and valuable suggestions on the manuscript.
They also acknowledge fruitful and interesting discussions with Claudia T. Code\c{c}o, Claudio J. Struchiner and Daniel A. M. Villela.

\bibliographystyle{siam}
\bibliography{biblio}

\end{document}